\documentclass[12inch,a4paper]{article}
\usepackage[top=3cm, bottom=2.5cm, left=2.4cm, right=2.7cm]{geometry}
\usepackage{graphicx}
\usepackage{epstopdf}
\usepackage{booktabs}
\usepackage{blkarray}
\usepackage[british]{babel}

 \usepackage[sort,comma,authoryear]{natbib} 

    \setcitestyle{citesep={;}}
\setcitestyle{square}
\usepackage{amssymb}
\usepackage{setspace}
\usepackage{amsfonts}
\usepackage{amsmath}
\usepackage{float}
\usepackage{color}
\usepackage[dvipsnames]{xcolor}
 \usepackage{colortbl}
\usepackage{url}

\definecolor{Crimson}{rgb}{0.6471, 0.1098, 0.1882}

\usepackage{epstopdf}
\newtheorem{lemma}{Lemma}
\newtheorem{remark}{\normalfont \textit{Remark}}
\newtheorem{theorem}{Theorem}
\newtheorem{proposition}{Proposition}
\newtheorem{corollary}{Corollary}

\usepackage{tcolorbox}
\usepackage{epstopdf}
\usepackage{subcaption}
\providecommand{\keywords}[1]{\textit{Key words:} #1}
\usepackage[labelfont=bf]{caption}

\usepackage{titlesec}
\titleformat*{\subsection}{\normalfont \large}
\usepackage{abstract}
\usepackage{slashbox}
\begin{document}

% ********************************

% ********************************

\title{Consistency of the plug-in functional predictor of the Ornstein-Uhlenbeck process in Hilbert and Banach spaces}

\author{Javier \'Alvarez-Li\'ebana$^1$, Denis Bosq$^2$ and M. Dolores Ruiz--Medina$^1$}
\maketitle
\begin{flushleft}
$^1$ Department of Statistics and O. R., University of Granada, Spain. 
$^2$ LSTA, Universit\'e Pierre et Marie Curie--Paris 6, Paris, France.

\textit{E-mail: javialvaliebana@ugr.es}
\end{flushleft}

\doublespacing

% ********************************

% ********************************

\renewcommand{\absnamepos}{flushleft}
\setlength{\absleftindent}{0pt}
\setlength{\absrightindent}{0pt}
\renewcommand{\abstractname}{Summary}
\begin{abstract}
New results on functional prediction of the Ornstein-Uhlenbeck process in an autoregressive Hilbert-valued
and Banach-valued frameworks are derived. Specifically, consistency of the maximum likelihood estimator of the
autocorrelation operator, and of the associated plug-in predictor is obtained in both frameworks.

\vspace{0.5cm}
\textbf{Published in Statistics \& Probability Letters 117:12--22.}
 DOI: \url{doi.org/10.1016/j.spl.2016.04.023}
\end{abstract}


\textcolor{Crimson}{\section{Introduction}
\label{A2:sec:1}}

This paper derives new   results in the context of linear processes in function spaces.  An extensive literature has been developed in this context in the last few decades (see, for example, \cite{Bosq00,FerratyVieu06,RamsaySilverman05};
among others). In particular, the problem of functional prediction of linear processes  in Hilbert and Banach spaces has been widely addressed. We refer to the reader to the papers by \cite{BensmainMourid01}, \cite{Bosq96,Bosq02,Bosq04,Bosq07}, \cite{Guillas00,Guillas01}, \cite{Mas02,Mas04,Mas07}, \cite{MasMennetau03a,Menneteau05}, \cite{LabbasMourid02,MokhtariMourid03,Mourid02,Mourid04}
\cite{Rachedi04,Rachedi05,RachediMourid03}, \cite{DedeckerMerlevede03,DehlingSharipov05,GlendinningFleet07,Kargin08,Ruiz12}, \cite{Pumo98,MarionPumo04} and  \cite{Turbillonetal07,Turbillonetal08}; and the references
therein.  In the above--mentioned papers, different projection
methodologies have been adopted in the derivation of the main
asymptotic properties of the formulated  functional parameter
estimators and predictors. Particularly, \cite{Bosq00,BosqBlanke07} apply
Functional Principal Component Analysis (FPCA); \cite{AntoniadisSapatinas03,Antoniadisetal06,LaukaitisVasilecas09} propose wavelet--bases--based estimation methods. Applications  of these functional estimation results can be found in
the papers by \cite{DamonGuillas02,AntoniadisSapatinas03,Laukaitis08,HormannKokoszka11,RuizSalmeron09}; among others.

 We here pay attention  to the problem of
functional prediction of the Ornstein--Uhlenbeck (O.U.) process (see,
for example, \cite{UhlenbeckOrnstein30,WangUhlenbeck45}, for its introduction
and properties). See also \cite{Doob42} for the
classical definition of O.U. process from the Langevin (linear)
stochastic differential equation. We can find in
\cite{Kutoyants04,LiptserShiraev01} an explicit
expression of the maximum likelihood estimator (MLE) of the scale
parameter $\theta ,$ characterizing its covariance function. Its
strong consistency is proved, for instance, in \cite{KleptsynaBreton02}. We formulate  here the O.U. process as an  autoregressive Hilbertian process of
order one (so--called ARH(1) process),  and as an autoregressive Banach--valued
process of order one (so--called ARB(1) process). Consistency of the MLE of
$\theta $  is applied to prove the consistency of the corresponding MLE
of the autocorrelation operator of the O.U. process. We adopt the
methodology applied in \cite{Bosq91}, since our
interest relies on forecasting the values of the O.U. process over
an entire time interval. Specifically, considering  the O.U. process
$\left\lbrace \xi_{t}, \ t\in \mathbb{R} \right\rbrace$ on the basic
probability space $(\Omega, \mathcal{A},\mathcal{P}),$ we can define

\begin{equation}
X_{n}(t)=\xi_{nh+t},\quad 0\leq t \leq h,\quad n\in \mathbb{Z},
\label{A2:eqarh1OU} 
\end{equation} 
\noindent satisfying

\begin{equation}
X_{n} \left(t\right) = \xi_{nh + t} = \int_{-\infty}^{nh+t}
e^{-\theta \left(nh+t - s \right)} dW_s = \rho_{\theta}
\left(X_{n-1}\right) \left(t\right)   + \varepsilon_{n}
\left(t\right), \quad n \in \mathbb{Z}, \label{A2:39} 
\end{equation}
\noindent with

\begin{eqnarray}
\rho_{\theta} \left(x\right)\left(t\right) &=& e^{-\theta t} x\left(h
\right), \quad \rho_{\theta} \left(X_{n-1} \right)\left(t\right) =
e^{-\theta t} \int_{-\infty}^{nh} e^{-\theta \left(nh - s \right)}
dW_s, \nonumber \\
 \varepsilon_{n} \left(t\right) &=& \int_{nh}^{nh+t}
e^{-\theta \left(nh + t - s \right)} dW_s, \nonumber \\ \label{A2:40}
\end{eqnarray}

\noindent for $0 \leq t \leq h$, where $W = \left \lbrace W_t, \ t \in \mathbb{R} \right\rbrace$ is a
standard bilateral Wiener process (see \textcolor{Crimson}{Supplementary Material} \ref{A2:Supp}). Thus, $X=\left\lbrace X_{n},\ n\in \mathbb{Z} \right\rbrace$
satisfies the ARH(1) equation (\ref{A2:39}) (see also equation
(\ref{A2:24}) below for its general definition). The real separable
Hilbert space $H$ is given by \linebreak $H =
L^2\left(\left[0,h\right],\beta_{\left[0,h\right]},\lambda +
\delta_{(h)} \right),$ where
 $\beta_{\left[0,h\right]}$ is the Borel $\sigma$-algebra generated
by the subintervals in $\left[0,h\right],$  $\lambda$ is the
Lebesgue measure and $\delta_{(h)}(s)= \delta \left(s - h \right)$
is the Dirac measure at point $h.$ The associated norm
$$\|f\|_{H}= \sqrt{\displaystyle \int_{0}^{h}\left(f(t)\right)^2dt+ \left(f(h) \right)^2}, \quad f \in H =  L^2\left(\left[0,h\right],\beta_{\left[0,h\right]}, \lambda +
\delta_{(h)}\right),$$
\noindent  establishes the equivalent classes of functions given by
the relationship $f \sim_{\lambda + \delta_{(h)}} g$ if and only if
 $$\left(\lambda +\delta_{(h)} \right) \left(\left\lbrace
t:f\left(t\right) \neq g\left(t\right)\right\rbrace\right)  = 0,$$
with
\begin{eqnarray}
\left(\lambda +\delta_{(h)} \right) \left(\left\lbrace
t:f\left(t\right) \neq g\left(t\right)\right\rbrace\right)  = 0
\Leftrightarrow \lambda\left(\left\lbrace t:f\left(t\right) \neq
g\left(t\right)\right\rbrace\right) = 0 \text{ and } f\left(h\right)
= g\left(h\right), \nonumber %\label{A2:41}
\end{eqnarray}
\noindent where, as before,  $\delta_{(h)}$ is the Dirac measure at point $h$. We
will prove, in \textcolor{Crimson}{Lemma} \ref{A2:lem1} below,  that \linebreak $X= \left\lbrace X_{n},\ n\in \mathbb{Z} \right\rbrace,$ constructed in (\ref{A2:eqarh1OU}) from the O.U. process,
satisfying equations (\ref{A2:39})--(\ref{A2:40}), is  the unique
stationary solution to equation (\ref{A2:39}), in the space $H=L^2
\left(\left[0,h\right],\beta_{\left[0,h\right]},\lambda +
\delta_{(h)}\right),$ admitting a MAH($\infty$) representation.
Similarly,  in \textcolor{Crimson}{Lemma} \ref{A2:lem2} below, we will prove that
$X= \left\lbrace X_{n},\ n\in \mathbb{Z} \right\rbrace$, constructed in (\ref{A2:eqarh1OU}) from
the O.U. process, satisfying equations (\ref{A2:39})--(\ref{A2:40}), is
the unique stationary solution to equation (\ref{A2:39}),  admitting a
MAB($\infty$) representation, in the space $B =
\mathcal{C}\left(\left[0,h\right] \right),$  the real separable Banach space of
continuous functions, whose support is the interval
$\left[0,h\right],$ with the supremum norm.

The main results of this paper provide the almost surely convergence
to $\rho_{\theta }$ of its MLE $\rho_{\widehat{\theta }}$,
 in the norm of $\mathcal{L}(H),$ the space of bounded linear operators in the Hilbert space $H$ (respectively,
in the norm of $\mathcal{L}(B),$ the space of bounded linear
operators in the Banach space $B$). The convergence in probability
of the associated plug--in ARH(1) and ARB(1) predictors (i.e., the
convergence in probability of $\rho_{\widehat{\theta }}(X_{n-1})$ to
$\rho_{\theta }(X_{n-1})$ in $H$ and $B,$ respectively) is proved as
well.

 The
outline of this paper is as follows.  In \textcolor{Crimson}{Appendix} \ref{A2:sec:5}, the
main results of this paper are obtained. Specifically, \textcolor{Crimson}{Appendix} \ref{A2:s21} provides the definition of an O.U. process as an ARH(1)
process. Strong consistency in $\mathcal{L}(H)$ of the  estimator of
the autocorrelation operator is derived in \textcolor{Crimson}{Appendix} \ref{A2:s22}. Consistency in $H$ of the associated  plug--in ARH(1) predictor is
then established in \textcolor{Crimson}{Appendix} \ref{A2:s23}. The corresponding results in
Banach  spaces are given in \textcolor{Crimson}{Appendix} \ref{A2:sec:52}. For illustration
purposes, a simulation study is undertaken in \textcolor{Crimson}{Appendix} \ref{A2:sec:9}.
Final comments can be found in \textcolor{Crimson}{Appendix} \ref{A2:sec:7}. The basic
preliminary elements, applied in the proof of the main results of
this paper, and the proof of \textcolor{Crimson}{Lemma} \ref{A2:lem1},  can be found in the
\textcolor{Crimson}{Supplementary Material} \ref{A2:Supp}.

\textcolor{Crimson}{\section{Prediction of O.U. processes in Hilbert and Banach spaces}
\label{A2:sec:5}}

 In this section, we consider $H$ to be a real
separable Hilbert space. Recall that  a zero--mean  ARH(1) process
$X= \left\lbrace X_{n},\ n\in \mathbb{Z} \right\rbrace$, on the basic probability space
$(\Omega, \mathcal{A},\mathcal{P}),$  satisfies (see \cite{Bosq00})
\begin{equation}
X_{n } (t) = \rho \left(X_{n-1} \right) (t) + \varepsilon_{n}(t),\quad n \in
\mathbb{Z}, \quad \rho \in \mathcal{L}(H),\label{A2:24}
\end{equation}
\noindent  where $\rho $ denotes the autocorrelation operator of
process $X.$ Here, $\varepsilon= \left\lbrace \varepsilon_{n }, \ n\in
\mathbb{Z} \right\rbrace$ is assumed to be a strong--white noise; i.e.,
$\varepsilon $ is a Hilbert--valued zero-mean stationary process,
with independent and identically  distributed components in time,
 with $\sigma^{2}={\rm E} \left\lbrace \|\varepsilon_{n}\|_{H}^{2} \right\rbrace<\infty,$ for all
$n\in \mathbb{Z}.$

\textcolor{Crimson}{\subsection{O.U. processes as ARH(1) processes}
\label{A2:s21}}

 As commented in \textcolor{Crimson}{Appendix} \ref{A2:sec:1}, equations
(\ref{A2:eqarh1OU})--(\ref{A2:40}) provide the definition of an O.U. process
as an ARH(1) process, with $H =
L^2\left(\left[0,h\right],\beta_{\left[0,h\right]},\lambda +
\delta_{(h)} \right).$    The norm in the space $H$ of $\rho_{\theta }(x),$    with $\rho_{\theta
}$ introduced in (\ref{A2:40}) and  $x\in H,$   is given by
\begin{equation}
\Vert \rho_{\theta}(x) \Vert_{H}^{2}  = \displaystyle \int_{0}^{h}
\left( \rho_{\theta} \left(x\right)\left(t\right) \right)^2 d
\left(\lambda + \delta_{(h)}\right)\left(t\right) = \displaystyle
\int_{0}^{h} \left( \rho_{\theta} \left(x\right)\left(t\right)
\right)^2 dt + \left(\rho_{\theta} \left(x\right)\left(h\right)
\right)^2, \nonumber %\label{A2:42}
\end{equation}
\noindent for each $h > 0$. The following lemma
provides,  for each $k\geq 1,$  the exact value of the norm of
$\rho^{k}_{\theta },$ in the space of bounded linear operators on
$H.$  As a direct consequence, the existence of an integer
$k_{0}$ such that $\Vert \rho_{\theta}^{k}
\Vert_{\mathcal{L}\left(H\right)}<1,$  for $k\geq k_{0},$ is also
derived for $\theta
>0.$

\bigskip

\begin{lemma}
\label{A2:lem1} \textit{Let us consider $\theta > 0$ and $X = \left\lbrace X_n, \ n\in \mathbb{Z} \right\rbrace$
satisfying equations (\ref{A2:eqarh1OU})--(\ref{A2:40}). For each $k\geq 1,$ the uniform norm of
$\rho_{\theta}^{k}$ is given by
\begin{equation}
\Vert\rho_{\theta}^{k}\Vert_{\mathcal{L}\left(H\right)}=\sqrt{e^{-2\theta
(k-1)h} \left(\frac{1+e^{-2\theta  h}\left(2\theta
-1\right)}{2\theta}\right)}=e^{-\theta (k-1)h}\Vert \rho_{\theta}
\Vert_{\mathcal{L}\left(H\right)}.\label{A2:nrkr}
\end{equation}
\noindent Furthermore,  for $k \geq
k_0=\left[\frac{1}{\theta} + 1\right]^{+},$
\begin{equation}
\Vert \rho_{\theta}^{k} \Vert_{\mathcal{L}\left(H\right)} < 1,\label{A2:eqL1}
\end{equation}
\noindent where $[t]^{+}$ denotes the closest upper integer of $t,$ for every $t\in \mathbb{R}_{+}.$}
\end{lemma}

\bigskip

The proof of this lemma can be found in the \textcolor{Crimson}{Supplementary Material} \ref{A2:Supp} provided.

\bigskip

\begin{remark}
\label{A2:R1} \textit{From equation (\ref{A2:eqL1}), applying 
\cite[Theorem 3.1]{Bosq00},  \textcolor{Crimson}{Lemma} \ref{A2:lem1} implies
that $X$  constructed in (\ref{A2:eqarh1OU})
from an O.U. process, defines the unique stationary solution to
equation (\ref{A2:39}) in the space  \linebreak $H= L^2
\left(\left[0,h\right],\beta_{\left[0,h\right]},\lambda +
\delta_{(h)}\right),$ admitting the MAH($\infty $) representation
\begin{equation}
X_{n } = \displaystyle\sum_{k=0}^{+ \infty} \rho_{\theta}^{k}
\left(\varepsilon_{n -k} \right), \quad n \in \mathbb{Z}, \quad\rho_{\theta} \in
\mathcal{L} \left(H\right). \nonumber %\label{A2:25}
\end{equation}}
\end{remark}

\bigskip

\begin{remark}
\textit{Note  that, for all $x\in H,$ and $k\geq 2,$ $\Vert
\rho_{\theta}^{k}\Vert_{\mathcal{L}\left(H\right)}\leq \Vert
\rho_{\theta}\Vert_{\mathcal{L}\left(H\right)}^{k}.$}
\end{remark}

\textcolor{Crimson}{\subsection{Functional parameter estimation and consistency}
\label{A2:s22}}

We now prove the strong consistency of the estimator
$\rho_{\widehat{\theta}_{n}}$  of operator $\rho_\theta$ in
$\mathcal{L}(H),$ with, as before, \linebreak $H= L^2
\left(\left[0,h\right],\beta_{\left[0,h\right]},\lambda +
\delta_{(h)}\right),$ and  $\widehat{\theta}_{n}$ denoting  the  MLE
of $\theta,$ based on the observation of an O.U. process on the
interval $[0,T],$ with $T=nh.$ Note that, from equation (\ref{A2:40}),
for all $x\in H,$ and for a given sample size $n,$
$$\rho_{\widehat{\theta}_{n}}(x)= e^{-\widehat{\theta}_{n} t}
x\left(h \right),$$ where  the  MLE of $\theta $ is given, for
$T=nh,$ by
\begin{eqnarray} \widehat{\theta}_{T} &=& \frac{1 +
\frac{\xi_{0}^{2}}{T} -
\frac{\xi_{T}^{2}}{T}}{\frac{2}{T}\displaystyle \int_{0}^{T}\xi_{t}^{2} dt}, \quad
T > 0,\label{A2:17b}
\end{eqnarray}
\noindent with $\left\lbrace \xi_{t}, \ t\in [0,T] \right\rbrace$ being the observed values of
the O.U. process  over the interval $[0,T].$ Thus,
$\rho_{\widehat{\theta}_{n}}$ is introduced in an abstract way,
since it can only be explicitly computed, for each particular
function $x\in H$ considered.  However, the norm $\Vert \rho_\theta
- \rho_{\widehat{\theta}_n}\Vert_{\mathcal{L} \left(H \right)}$ is
explicitly
 computed in equation (\ref{A2:48}) below.

 The following results will be applied in the proof of \textcolor{Crimson}{Proposition} \ref{A2:prf}.
 
 \bigskip
 
\begin{lemma}
\label{A2:lem000} \textit{If $t \in \left[0, + \infty \right)$, it holds that
$$\vert e^{- u t} - e^{-v t} \vert \leq \vert u - v \vert t, \quad u,v \geq 0.$$}
\end{lemma}

 \bigskip

The proof of this lemma is given  in the \textcolor{Crimson}{Supplementary Material} \ref{A2:Supp}.

 \bigskip

\begin{theorem} \label{A2:thpr1} 
\textit{(See  also \cite[Proposition 2.2]{KleptsynaBreton02} and
\cite[p. 63 and p. 117]{Kutoyants04}).  The MLE of $\theta$ defined
in equation (\ref{A2:17b}) is strongly consistent; i.e.,
\begin{equation}
\widehat{\theta}_T \longrightarrow \theta \quad a.s., \quad T \rightarrow \infty. \nonumber %\label{A2:19a}
\end{equation}}
\end{theorem}

 \bigskip
The proof follows from the Ibragimov--Khasminskii's Theorem.
\bigskip

\begin{proposition}
\label{A2:prf} 
\textit{Let $H$ be the space $L^2\left(\left[0,h\right],
\beta_{\left[0,h\right]}, \lambda + \delta_{(h)}\right).$  Then, the
estimator $\rho_{\widehat{\theta}_{n}}$ of operator $\rho_{\theta},$
based on the MLE $\widehat{\theta}_{n}$ of $\theta$, is strongly
consistent in the norm of $\mathcal{L} \left(H \right)$; i.e.,
\begin{equation}
\Vert \rho_\theta - \rho_{\widehat{\theta}_n} \Vert_{\mathcal{L}
\left(H \right)} \longrightarrow 0 \quad a.s., \quad n \rightarrow \infty. \nonumber %\label{A2:47}
\end{equation}}
\end{proposition}

\begin{proof}
The following straightforward almost surely identities are obtained:
\begin{eqnarray}
\Vert \rho_\theta - \rho_{\widehat{\theta}_n} \Vert_{\mathcal{L}
\left(H \right)}&=& \displaystyle \sup_{x\in H} \left\lbrace
\frac{\Vert \left( \rho_{\theta} - \rho_{\widehat{\theta}_n}
\right)\left(x \right)
\Vert_H}{\Vert x  \Vert_H} \right\rbrace \nonumber \\
&=& \displaystyle \sup_{x\in H}
\left\lbrace\sqrt{\frac{\displaystyle \int_{0}^{h} \left(
\left(\rho_{\theta} - \rho_{\widehat{\theta}_n} \right)\left(x
\right) \left( t \right)\right)^2 d\left(\lambda
+ \delta_{(h)}\right) \left( t \right)}{\displaystyle \int_{0}^{h} \left( x\left(t \right)\right)^2 d\left(\lambda + \delta_{(h)}\right) \left( t \right)}}\right\rbrace \nonumber\\
&=&  \displaystyle \sup_{x\in H} \left\lbrace\sqrt{ \left(x \left(h
\right) \right)^2 \frac{\displaystyle \int_{0}^{h} \left(e^{-\theta t} -
e^{-\widehat{\theta}_n t} \right)^2 dt + \left(e^{-\theta h} -
e^{-\widehat{\theta}_n h} \right)^2} {\displaystyle \int_{0}^{h} \left(x \left(t
\right) \right)^2 dt + \left(x \left(h
\right) \right)^2}}\right\rbrace \nonumber
\end{eqnarray}
\begin{eqnarray}
&=& \sqrt{\displaystyle \int_{0}^{h} \left(e^{-\theta t} -
e^{-\widehat{\theta}_n t} \right)^2 dt + \left(e^{-\theta h} -
e^{-\widehat{\theta}_n h} \right)^2}, \label{A2:48}
\end{eqnarray}
\noindent where the last identity is obtained in a similar way to
equation (\ref{A2:nrkr}) in \textcolor{Crimson}{Lemma}
\ref{A2:lem1} (see \textcolor{Crimson}{Supplementary Material} \ref{A2:Supp}).

 From \textcolor{Crimson}{Lemma} \ref{A2:lem000}  and equation (\ref{A2:48}), for $n$ sufficiently large, we have
\begin{eqnarray}
\Vert \rho_\theta - \rho_{\widehat{\theta}_n}\Vert_{\mathcal{L}
\left(H \right)} &\leq & \sqrt{\displaystyle \int_{0}^{h} t^{2}\vert
\theta -\widehat{\theta}_{n}\vert^{2}dt +
h^{2}\vert \theta -\widehat{\theta}_{n}\vert^{2}} = \vert \theta -\widehat{\theta}_{n}\vert\sqrt{\displaystyle \int_{0}^{h} t^{2}dt + h^{2}}  \nonumber \\
&=& \vert \theta -\widehat{\theta}_{n}\vert h \sqrt{\frac{h}{3} +
1} \quad a.s. \label{A2:49}
\end{eqnarray}

The strong--consistency of $\rho_{\widehat{\theta}_n}$ in
$\mathcal{L} \left(H \right)$ directly follows from \textcolor{Crimson}{Theorem} \ref{A2:thpr1} and  equation (\ref{A2:49}).

\hfill \hfill  \textcolor{Aquamarine}{$\blacksquare$}
\end{proof}

\bigskip

\begin{remark}
\label{A2:rem3b}\textit{ From 
\cite[Proposition 2.3]{KleptsynaBreton02} (see
also \textcolor{Crimson}{Theorem} \ref{A2:th3prem}  below), the MLE $\widehat{\theta}_T$ of
$\theta $ satisfies
\begin{equation} {\rm E}\left\lbrace \left( \theta - \widehat{\theta}_T \right)^2 \right\rbrace =\mathcal{O} \left( \frac{ 2 \theta } {T} \right),\quad
T\rightarrow \infty.\label{A2:rnctt}
\end{equation} \noindent  In addition, from equation (\ref{A2:49}), considering  $T=nh,$ $h>0,$
\begin{equation}{\rm E} \left\lbrace \Vert \rho_\theta -
\rho_{\widehat{\theta}_n}\Vert^{2}_{\mathcal{L} \left(H
\right)}\right\rbrace \leq {\rm E} \left\lbrace \vert \theta
-\widehat{\theta}_{n}\vert^{2}\right\rbrace h^{2} \left(\frac{h}{3} +
1\right).\label{A2:rncb}
\end{equation}
\noindent Equations (\ref{A2:rnctt})--(\ref{A2:rncb}) lead to
$${\rm E} \left\lbrace \Vert \rho_\theta -
\rho_{\widehat{\theta}_n}\Vert^{2}_{\mathcal{L} \left(H
\right)}\right\rbrace \leq G(\theta, \widehat{\theta }_{n},h),$$ \noindent
with $$G(\theta, \widehat{\theta
}_{n},h)=\mathcal{O}\left(\frac{2\theta}{n}\right), \quad n\rightarrow
\infty.$$
Therefore, the functional parameter  estimator
$\rho_{\widehat{\theta}_n}$ is $\sqrt{n}$--consistent.}
\end{remark}

\bigskip

\textcolor{Crimson}{\subsection{Consistency of the plug--in ARH(1) predictor}
\label{A2:s23}}
 Let us consider the plug--in ARH(1) predictor
$\widehat{X}_{n},$ constructed from the MLE
$\rho_{\widehat{\theta}_n}$ of $\rho_{\theta }$
 in \textcolor{Crimson}{Proposition} \ref{A2:prf}, given by
\begin{equation}
\widehat{X}_{n} \left( t \right)= \rho_{\widehat{\theta}_n}
\left(X_{n-1}\right)\left(t\right) = e^{-\widehat{\theta}_n t}
X_{n-1} \left( h \right),\quad 0 \leq t \leq h, \quad n \in \mathbb{Z}. \label{A2:50}
\end{equation}

\textcolor{Crimson}{Corollary} \ref{A2:c1} below provides the consistency of
$\widehat{X}_{n},$ given in equation (\ref{A2:50}), from  \textcolor{Crimson}{Proposition}
\ref{A2:prf}  by applying the following lemma and theorem.

\bigskip

\begin{lemma} 
\label{A2:lem2pre} 
\textit{Let $\left\lbrace Z_n, \ n \in \mathbb{Z} \right\rbrace$ be a
sequence of random variables such that $$Z_n \sim \mathcal{N}
\left(0, \frac{1}{2 \theta} \right), \quad \theta > 0,$$ and let
$\left\lbrace Y_n, \ n \in \mathbb{Z} \right\rbrace$ be another  sequence of
random variables such that $$\sqrt{\ln \left( n \right) }Y_n \longrightarrow^{p}
0, \quad n\rightarrow \infty.$$ Then, $$Y_n \vert Z_n \vert
\longrightarrow^{p} 0, \quad n\rightarrow \infty,$$ where, as usual, $\longrightarrow^{p}$
indicates convergence in probability.}
\end{lemma}

\bigskip
The proof of this lemma can be found in the \textcolor{Crimson}{Supplementary Material} \ref{A2:Supp}.

\bigskip

\begin{theorem}
\label{A2:th3prem}  
\textit{Let $\widehat{\theta}_T$ be the MLE of $\theta$
defined in equation (\ref{A2:17b}), with $\theta > 0.$ Hence,
\begin{equation}
 {\rm E} \left\lbrace \left( \theta - \widehat{\theta}_T \right)^2 \right\rbrace =\mathcal{O} \left( \frac{ 2 \theta } {T} \right),\quad T\rightarrow \infty .\label{A2:18a}
\end{equation}
In particular,  \begin{equation}
\lim_{T\rightarrow \infty } {\rm E} \left\lbrace
\left( \theta - \widehat{\theta}_T \right)^2
\right\rbrace=0. \nonumber %\label{A2:eqaemlee}
\end{equation}}
\end{theorem}

\bigskip

The proof of this result is given in  \cite[Proposition 2.3]{KleptsynaBreton02}.

\bigskip

\begin{corollary}
\label{A2:c1}
\textit{Let $H =
L^2\left(\left[0,h\right],\beta_{\left[0,h\right]},\lambda +
\delta_{(h)} \right)$ be the Hilbert space introduced above. Then,
the   plug--in ARH(1) predictor (\ref{A2:50})   of an O.U. process is
consistent in $H$; i.e.,
\begin{equation}
\left\| \left(\rho_\theta - \rho_{\widehat{\theta}_n} \right)
\left(X_{n-1} \right) \right\|_H \longrightarrow^{p} 0. \nonumber %\label{A2:51}
\end{equation}}
\end{corollary}

\begin{proof}
By definition,
\begin{eqnarray}
\left\| \left( \rho_\theta - \rho_{\widehat{\theta}_n}
\right)\left(X_{n-1} \right) \right\|_{H}  = \left| X_{n-1}\left(h
\right) \right| \sqrt{  \displaystyle \int_{0}^{h} \left(e^{-\theta t}
-  e^{-\widehat{\theta}_n t}\right)^2 dt + \left(e^{-\theta h} -
e^{-\widehat{\theta}_n h}\right)^2 }. \label{A2:52}
\end{eqnarray}
From equations (\ref{A2:48})--(\ref{A2:49}) and (\ref{A2:52}), we then
obtain, for $n$ sufficiently large,
\begin{eqnarray}
\left\|  \left( \rho_\theta - \rho_{\widehat{\theta}_n}
\right)\left(X_{n-1 } \right) \right\|_{H} \leq   \left| X_{n-1}\left( h
\right) \right| \left| \theta -\widehat{\theta}_n \right| h
\sqrt{\frac{h}{3} + 1}\quad a.s. \label{A2:53}
\end{eqnarray}

Let us set $$\left\lbrace Y_n, \ n \in \mathbb{Z} \right\rbrace = \left\lbrace
\vert \theta -\widehat{\theta}_n \vert h \sqrt{\frac{h}{3} + 1}, \ n \in \mathbb{Z}\right\rbrace, \quad \left\lbrace Z_n, \ n \in
\mathbb{Z}\right\rbrace = \left\lbrace   X_{n-1} \left( h \right), \ n \in
\mathbb{Z}\right\rbrace,$$ with $Z_n \sim \mathcal{N} \left(0, \frac{1}{2
\theta} \right),$ for every $n\in \mathbb{Z}.$ From \textcolor{Crimson}{Theorem}
\ref{A2:thpr1}, $$Y_n \longrightarrow 0\quad a.s., \quad n \rightarrow \infty .$$ Hence, to
apply \textcolor{Crimson}{Lemma} \ref{A2:lem2pre}, we need to prove that
\begin{equation}
\sqrt{\ln \left(n \right)} Y_n \longrightarrow^{p} 0,\quad n\rightarrow \infty. \nonumber %\label{A2:54}
\end{equation}

From the Chebyshev's inequality and \textcolor{Crimson}{Theorem} \ref{A2:th3prem}, we get, for
all $\varepsilon > 0,$
\begin{equation}
\displaystyle \lim_{n \rightarrow 0} \mathcal{P} \left( \vert \theta -\widehat{\theta}_n \vert \sqrt{\ln
\left(n \right)} h \sqrt{\frac{h}{3} + 1} \geq \varepsilon \right)
\leq  \frac{h^2 \left(\frac{h}{3} + 1\right) \ln\left( n \right)
{\rm E} \left\lbrace \left| \theta -\widehat{\theta}_n \right|^2
\right\rbrace}{\varepsilon^2} = 0. \nonumber %\label{A2:55}
\end{equation}

Therefore, from \textcolor{Crimson}{Lemma} \ref{A2:lem2pre},  we obtain the convergence in
probability of $\left\| \left(\rho_\theta - \rho_{\widehat{\theta}_n}
\right) \left(X_{n-1} \right) \right\|_{H}$ to zero.

\hfill \hfill \textcolor{Aquamarine}{$\blacksquare$}
\end{proof}

\textcolor{Crimson}{\subsection{Prediction of O.U. processes in $B = \mathcal{C}\left(\left[0,h\right] \right)$}
\label{A2:sec:52}}

As before, let $B$ be now the Banach space of continuous functions,
whose support is the interval $\left[0,h\right]$, with the supremum
norm, denoted as $\mathcal{C}\left(\left[0,h\right]\right).$ The
following lemma states that $\|\rho_{\theta }^{k}\|_{\mathcal{L}
\left(B \right)}\leq 1,$ for $\theta >0,$ and for every $k\geq 1,$
with $\mathcal{L} \left(B \right)$ being the space of bounded linear
operators on the Banach space \linebreak $B=
\mathcal{C}\left(\left[0,h\right]\right),$ and $\rho_{\theta }$
being introduced in equation (\ref{A2:40}). Consequently, from \cite[Theorem 6.1]{Bosq00}, $X= \left\lbrace X_{n},\ n\in
\mathbb{Z} \right\rbrace,$ constructed in (\ref{A2:eqarh1OU}) from the O.U. process,
 defines the unique stationary
solution to equation (\ref{A2:39}), in the Banach space $B=
\mathcal{C}\left(\left[0,h\right]\right),$ admitting a MAB($\infty$)
representation.

\bigskip

\begin{lemma}
\label{A2:lem2} 
\textit{Let $\rho_{\theta }$ introduced in (\ref{A2:40}), defined
on  $B= \mathcal{C}\left(\left[0,h\right]\right).$ Then, for  $k\geq
1,$ $\|\rho_{\theta }^{k}\|_{\mathcal{L} \left(B \right)}\leq 1,$
with $\theta
>0.$}
\end{lemma}

\begin{proof}

From $$\rho_{\theta}^{k}(x)(t)= e^{-\theta t}e^{-\theta
(k-1)h}x(h),$$ \noindent  for each  $k\geq 1$ and  $\theta
>0,$ we have
\begin{eqnarray}
\left\| \rho_{\theta}^{k}  \right\|_{\mathcal{L}\left(B\right)} &=&
\displaystyle \sup_{x\in B} \left\lbrace \frac{\Vert
\rho_{\theta}^{k} \left(x \right)  \Vert_B}{\Vert x  \Vert_B}
\right\rbrace  = \displaystyle \sup_{x\in B}  \left\lbrace  \frac{\displaystyle \sup_{0\leq
t\leq h}  \left\lbrace \left| e^{-\theta
t} e^{-\theta (k-1)h} x(h)\right| \right\rbrace}{\displaystyle \sup_{0\leq t\leq h} \left| x(t) \right|} \right\rbrace \nonumber\\
&=& \displaystyle \sup_{x\in B}  \left\lbrace  \frac{\left| x(h) \right| e^{-\theta
(k-1)h} \displaystyle \sup_{0\leq t\leq h}  e^{-\theta t} }{\displaystyle \sup_{0\leq t\leq
h} \left| x(t) \right|} \right\rbrace \leq \displaystyle \sup_{x\in B}  \left\lbrace \frac{\left| x(h) \right| \displaystyle \sup_{0\leq
t\leq h}  e^{-\theta t} }{\left| x(h) \right|} \right\rbrace
\nonumber\\
&=& \displaystyle \sup_{0\leq t\leq h}   e^{-\theta t}  = 1. \label{A2:56}
\end{eqnarray}
\hfill \hfill \textcolor{Aquamarine}{$\blacksquare $}
\end{proof}

\bigskip

We now check the strong consistency of the MLE
$\rho_{\widehat{\theta}_{n}}$ of $\rho_{\theta}$ in
$\mathcal{L}(B).$ From equation (\ref{A2:56}),
\begin{equation}
\Vert \rho_{\theta} - \rho_{\widehat{\theta}_n}\Vert_{\mathcal{L}
\left(B\right)} \leq \displaystyle \sup_{0 \leq t \leq h} \left\lbrace \left| e^{-\theta t} - e^{-\widehat{\theta}_n t}\right| \right\rbrace \quad a.s. \nonumber %\label{A2:59}
\end{equation}

 From  \textcolor{Crimson}{Lemma} \ref{A2:lem000}, for $n$ sufficiently large, we
then have
\begin{equation}\Vert \rho_{\theta} - \rho_{\widehat{\theta}_n}
\Vert_{\mathcal{L} \left(B \right)} \leq h \left| \theta -
\widehat{\theta}_n \right| \quad a.s.\label{A2:eqscbanach}
\end{equation}
\noindent \textcolor{Crimson}{Theorem} \ref{A2:thpr1} then leads to the desired result on
strong consistency of the estimator $\rho_{\widehat{\theta}_n}$ of
 $\rho_\theta$ in $\mathcal{L}(B).$ Furthermore, from \textcolor{Crimson}{Theorem}  \ref{A2:th3prem} , in a similar way to
 \textcolor{Crimson}{Remark} \ref{A2:rem3b}, the $\sqrt{n}$--consistency of
 $\rho_{\widehat{\theta}_n}$ in $\mathcal{L} \left(B \right)$ also
 follows from equations (\ref{A2:18a}) and (\ref{A2:eqscbanach}).

Similarly  to \textcolor{Crimson}{Corollary} \ref{A2:c1}, in the following result,
 the consistency, in the Banach space \linebreak $B=C([0,h]),$ of the   plug--in
predictor (\ref{A2:50}) is obtained.

\begin{corollary} 
\label{A2:c2}
\textit{ The  ARB(1) plug--in predictor (\ref{A2:50}) of
a zero--mean O.U. process is consistent in \linebreak $B=C([0,h])$; i.e., as
$n\rightarrow \infty,$
\begin{equation}
\left\| \left(\rho_\theta - \rho_{\widehat{\theta}_n} \right)
\left(X_{n-1} \right) \right\|_B \longrightarrow^{p} 0. \nonumber %\label{A2:60}
\end{equation}}
\end{corollary}

\begin{proof}
From \textcolor{Crimson}{Lemma}  \ref{A2:lem000}, for $n$ sufficiently large, and for each
$h>0,$
\begin{equation}
\Vert \left(\rho_\theta  - \rho_{\widehat{\theta}_n}
\right)\left(X_{n-1} \right) \Vert_B = \displaystyle \sup_{0 \leq
t \leq h} \left\lbrace \left| e^{-\theta t} - e^{-\widehat{\theta}_n t} \right| \left| X_{n-1} \left( h \right) \right|\right\rbrace \leq h \vert \theta
-\widehat{\theta}_n \vert \vert X_{n-1} \left( h \right) \vert \quad a.s. \label{A2:61}
\end{equation}

As derived in the proof of \textcolor{Crimson}{Corollary} \ref{A2:c1}, from \textcolor{Crimson}{Theorem} \ref{A2:th3prem},  the random sequence \linebreak $\left\lbrace Y_n, \ n \in \mathbb{Z} \right\rbrace = \left\lbrace h \vert \theta
-\widehat{\theta}_n \vert, \ n \in \mathbb{Z}\right\rbrace$ is such that
$$\sqrt{\ln \left( n \right)} Y_n \leq \sqrt{\frac{h}{3} + 1} \sqrt{\ln
\left( n \right)} Y_n \longrightarrow^{p}0, \quad n\rightarrow \infty.$$ Moreover,  $
\left\lbrace Z_n, \ n \in \mathbb{Z} \right\rbrace= \left\lbrace X_{n-1}
\left( h \right), \ n \in \mathbb{Z} \right\rbrace$ is such that
$Z_n \sim \mathcal{N}\left(0,\frac{1}{2\theta
}\right).$  Lemma \ref{A2:lem2pre}  then leads,  as $n\rightarrow
\infty,$ to the desired convergence result from equation (\ref{A2:61}):

\begin{equation}
\Vert \left(\rho_\theta - \rho_{\widehat{\theta}_n}\right)
\left(X_{n-1} \right) \Vert_B \leq  Y_n \vert Z_n \vert \longrightarrow^{p} 0. \nonumber %\label{A2:62}
\end{equation}
\hfill \hfill \textcolor{Aquamarine}{$\blacksquare $}
\end{proof}

\textcolor{Crimson}{\section{Simulations}
\label{A2:sec:9}}

In this section, a simulation study is undertaken to illustrate the
asymptotic  results presented in this paper about the MLE
$\widehat{\theta}_{n}$ of $\theta ,$ and the consistency of the ML
functional parameter estimators of the autocorrelation operator, and
the associated  plug--in predictors, in the ARH(1) and  ARB(1)
frameworks.

\textcolor{Crimson}{\subsection{Estimation of the scale parameter $\theta$}
\label{A2:sec:91}}

On the simulation of the sample--paths of an O.U. process, an extension of
the Euler's method, the so--called Euler--Murayama's method  (see
\cite{KloedenPlaten92}) is applied, from
the Langevin  stochastic differential equation satisfied by the O.U.
process $\left\lbrace \xi_t,\  t\in \left[0,T \right] \right\rbrace$
\begin{equation}
d\xi_t = -\theta \xi_t + dW_t,\quad \theta >0, \quad t \in \left[0,T
\right], \quad \xi_0 = x_0.\label{A2:170}
\end{equation}

\begin{figure}[H]
\centering
\includegraphics[width=0.8\textwidth]{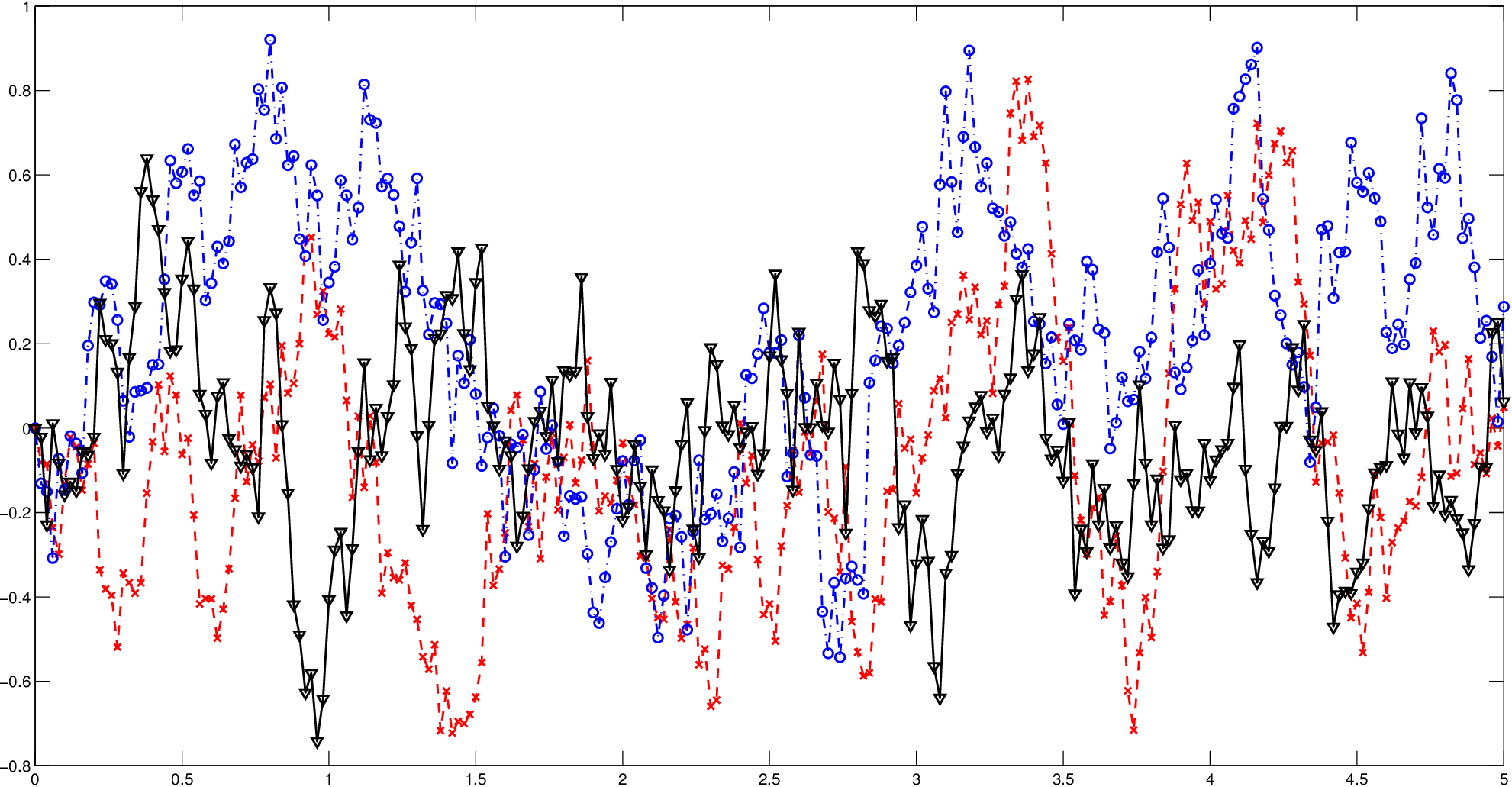}

\vspace{-0.1cm}
 \caption[\hspace{0.7cm} Sample paths of an O.U. process.]{{\small Sample paths of an O.U. process   $\left\lbrace \xi_t, \ 0 \leq t \leq T \right\rbrace$ generated with $T = 5$, $\Delta t = 0.02$, $\theta = 5$ and
    $\widehat{\xi}_0 = 0$.}}
  \label{A2:fig:Imagen1}
\end{figure}

Thus, let   $0=t_0 < t_1  < \dots < t_n = T$ be a partition of the
real interval $\left[0, T \right].$ Then, (\ref{A2:170}) can be
discretized as
\begin{equation}
\widehat{\xi}_{i+1} = \widehat{\xi}_{i} - \theta \widehat{\xi}_{i} +
\Delta W_i, \quad \widehat{\xi}_0 = \xi_0 = 0, \label{A2:171}
\end{equation}
where $\left\lbrace \Delta W_i, \ i=0,\dots,n-1 \right\rbrace$ are
i.i.d. Wiener increments; i.e., $$\Delta W_i \sim \mathcal{N}
\left(0,\Delta t \right) = \sqrt{\Delta t}  \mathcal{N} \left(0,1
\right), \quad  i=0,\dots,n-1.$$ In the following, we take $\Delta t = 0.02$ as
discretization step size, considering  $N = 1000$ simulations of the
 O.U. process. In particular, Figure \ref{A2:fig:Imagen1} shows some realizations
of the
 discrete version of the solution to (\ref{A2:170}) generated from  (\ref{A2:171}).

\begin{figure}[H]
\centering
    \includegraphics[width=0.8\textwidth]{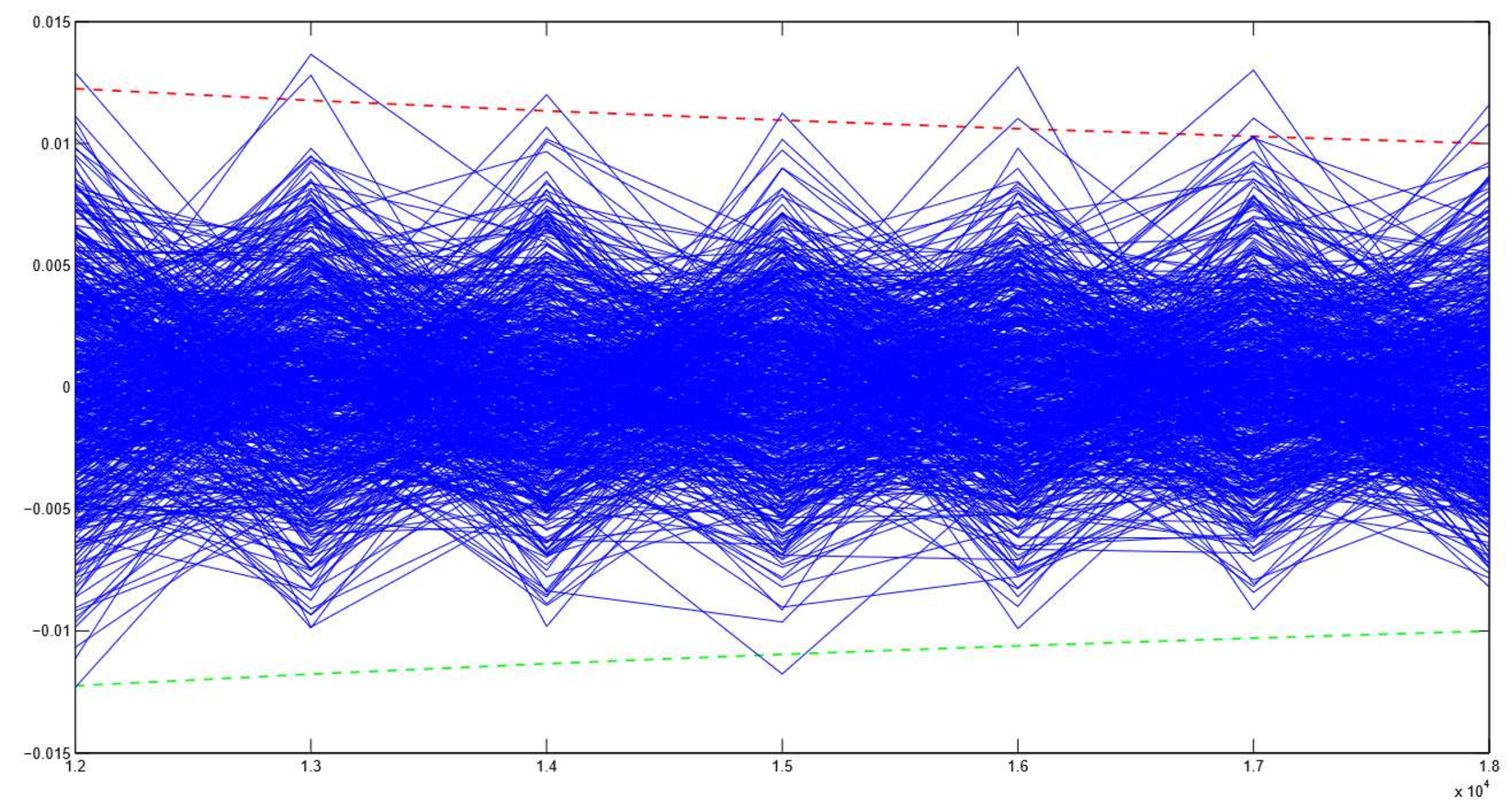}
  \label{A2:fig:Imagen3}

    \includegraphics[width=0.8\textwidth]{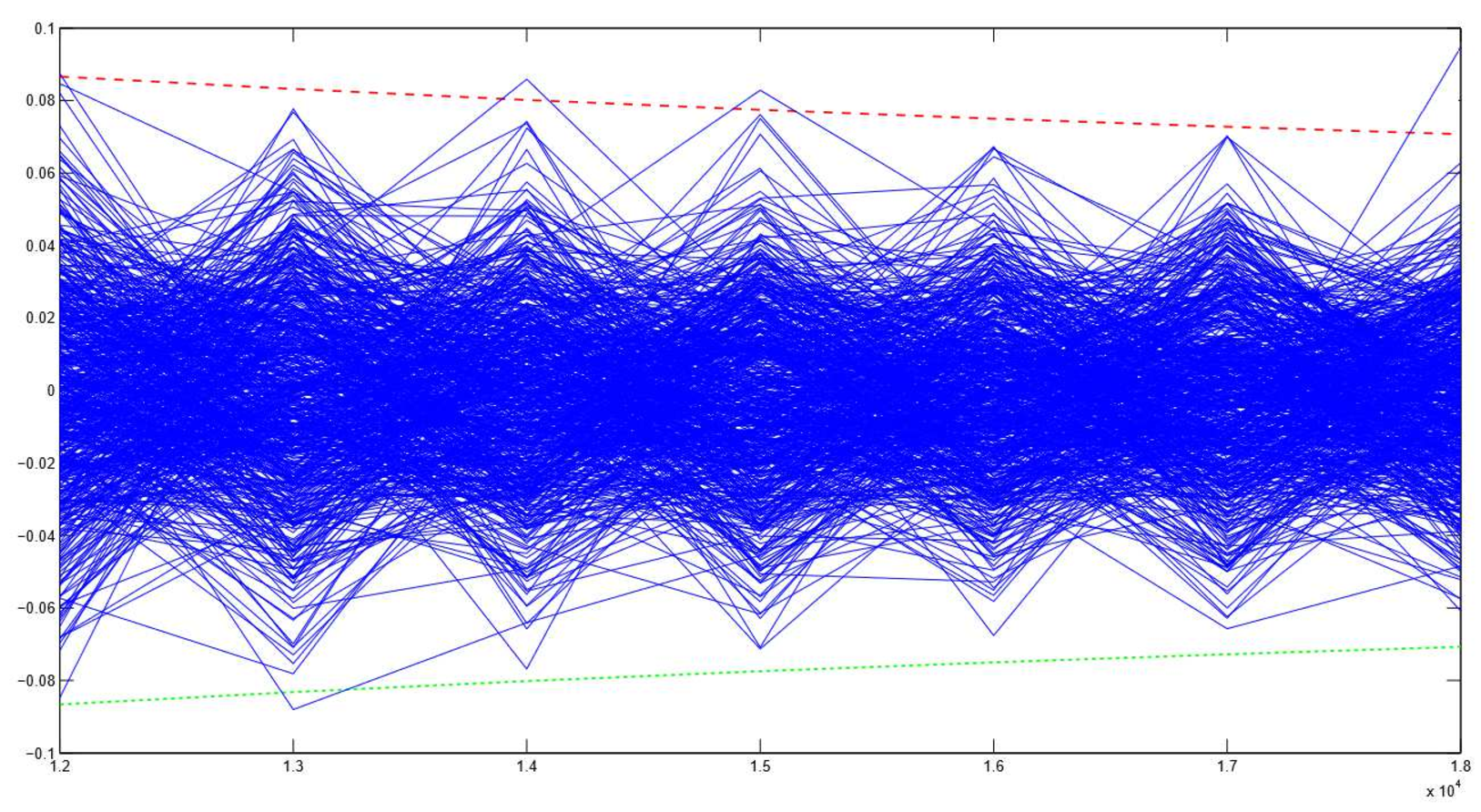}
  \label{A2:fig:Imagen4}
  
\vspace{-0.1cm}
\caption[\hspace{0.7cm} Empirical absolute errors on the estimation of $\theta$ of an O.U. process.]{{\small The values of $\widehat{\theta}_{T} - \theta ,$ based on $N = 1000$ simulations of the O.U. process over the interval $[0,T],$
for $ \left\lbrace T_l =12000+(l-1)1000, \ l=1,\dots,7 \right\rbrace,$ are represented  against
the confidence bands given by $+3\sigma=3\sqrt{\frac{2 \theta}{T}}$
(upper red dotted line) and $-3\sigma=-3\sqrt{\frac{2 \theta}{T}}$
(lower green dotted line), for  values $\theta= 0.1$ (at the top) and
$\theta= 5$ (at the bottom).}} \label{A2:fig:2}
\end{figure}

Let us first illustrate the asymptotic normal distribution of
$\widehat{\theta}_T$; i.e., for $T$ sufficiently large, we can
consider $\widehat{\theta}_T\sim \mathcal{N}
\left(\theta,\frac{2\theta }{ T}\right)$ (see \textcolor{Crimson}{Theorem} \ref{A2:Supp:th2prem} in the
\textcolor{Crimson}{Supplementary Material} \ref{A2:Supp}). From equation (\ref{A2:17b}),  we take
$$\widehat{\theta}_T = \frac{- \displaystyle \int_{0}^{T} \xi_t d
\xi_t}{\displaystyle \int_{0}^{T} \xi_{t}^{2} dt},$$ (see also
 \textcolor{Crimson}{Supplementary material} \ref{A2:Supp}), to compute the
following approximation of the MLE $\widehat{\theta}_{T}$ of $\theta
,$ for each one of the $N=1000$ simulations performed, and for each one of the six values of parameter
$\theta $ considered:

\begin{equation}
\widehat{\theta}_{T} \simeq \frac{- \displaystyle \sum_{i=0}^{n-1}
\widehat{\xi}_{t_{i},s}(\theta )
\left(\widehat{\xi}_{t_{i+1},s}(\theta ) -
\widehat{\xi}_{t_{i},s}(\theta ) \right)}{\displaystyle
\sum_{i=0}^{n-1} \widehat{\xi}_{t_{i},s}^{2}(\theta ) \Delta t}, \quad t_0
= 0,~t_n = T,~\Delta t = 0.02,\ ~s=1,\dots,N, \label{A2:173}
\end{equation}

\noindent where $\widehat{\xi}_{t_{i},s}(\theta )$ represents the
$s$--th discrete generation of the O.U. process,  evaluated at time
$t_i,$ with covariance scale parameter $\theta,$ for $$\theta = \left[ 0.1,
0.4, 0.7, 1, 2, 5\right].$$ Table \ref{A2:tab:Table1} displays the empirical
probabilities of the    error $\widehat{\theta}_{T} - \theta$ to be
within the band $\pm 3\sqrt{\frac{2 \theta}{T}},$ from    $N = 1000$
discrete simulations of  the O.U. process, considering different
sample sizes  \linebreak $\left\lbrace T_l = 12000+1000(l-1), \ l=1,\dots,7\right\rbrace$. Figure \ref{A2:fig:2} displays the
cases   $\theta =0.1$ (at the top) and $\theta =5$  (at
the bottom).
 It can be observed that, for each one of the sample sizes
 considered, $\left\lbrace T_l = 12000+1000(l-1), \ l=1,\dots,7\right\rbrace$,  approximately a 99\% of the realizations of  $\widehat{\theta}_{T} - \theta$ lie within the band $\pm 3\sqrt{\frac{2 \theta}{T}},$ which
supports the asymptotic Gaussian distribution.

\begin{table}[H]
\caption[\hspace{0.7cm} Empirical probabilities of the error of the MLE of the parameter of the O.U. process.]{\small{Empirical probabilities of the  error of the  MLE of
$\theta $ to lie within the band $\pm 3\sigma = \pm 3\sqrt{\frac{2
\theta}{T}},$
  for different sample
sizes $T,$  and values of parameter $\theta $.}}
\centering
\begin{small}
\begin{tabular}{|c||c|c|c|c|c|c|}
  \hline
   & \multicolumn{6}{c|}{Parameter $\theta$} \\
   \hline
   $T$ & $0.1$ & $0.4$  & $0.7$  & $1$ & $2$  & $5$   \\
  \hline \hline
   $12000$ & $0.998$ & $1$  & $0.998$ & $0.998$  & $1$ & $0.998$ \\
  \hline
   $13000$ & $0.997$ & $0.998$ & $0.998$ & $1$ & $0.995$ & $1$ \\
  \hline
     $14000$ & $0.998$ & $0.997$ & $1$ & $0.997$ & $1$ & $0.998$ \\
  \hline
     $15000$ & $0.998$ & $0.997$ & $0.998$ & $0.998$ & $1$ & $0.998$ \\
  \hline
     $16000$ & $0.997$ & $0.995$ & $0.997$ & $0.998$ & $1$ & $1$ \\
  \hline
    $17000$ & $0.993$ & $0.998$ & $1$ & $0.997$ & $0.995$ & $1$ \\
  \hline
     $18000$ & $0.997$ & $0.997$ & $0.995$ & $1$ & $1$ & $0.998$ \\
  \hline
\end{tabular}
\end{small}
  \label{A2:tab:Table1}
\end{table}

\begin{figure}[H]
   \centering
    \includegraphics[width=0.8\textwidth]{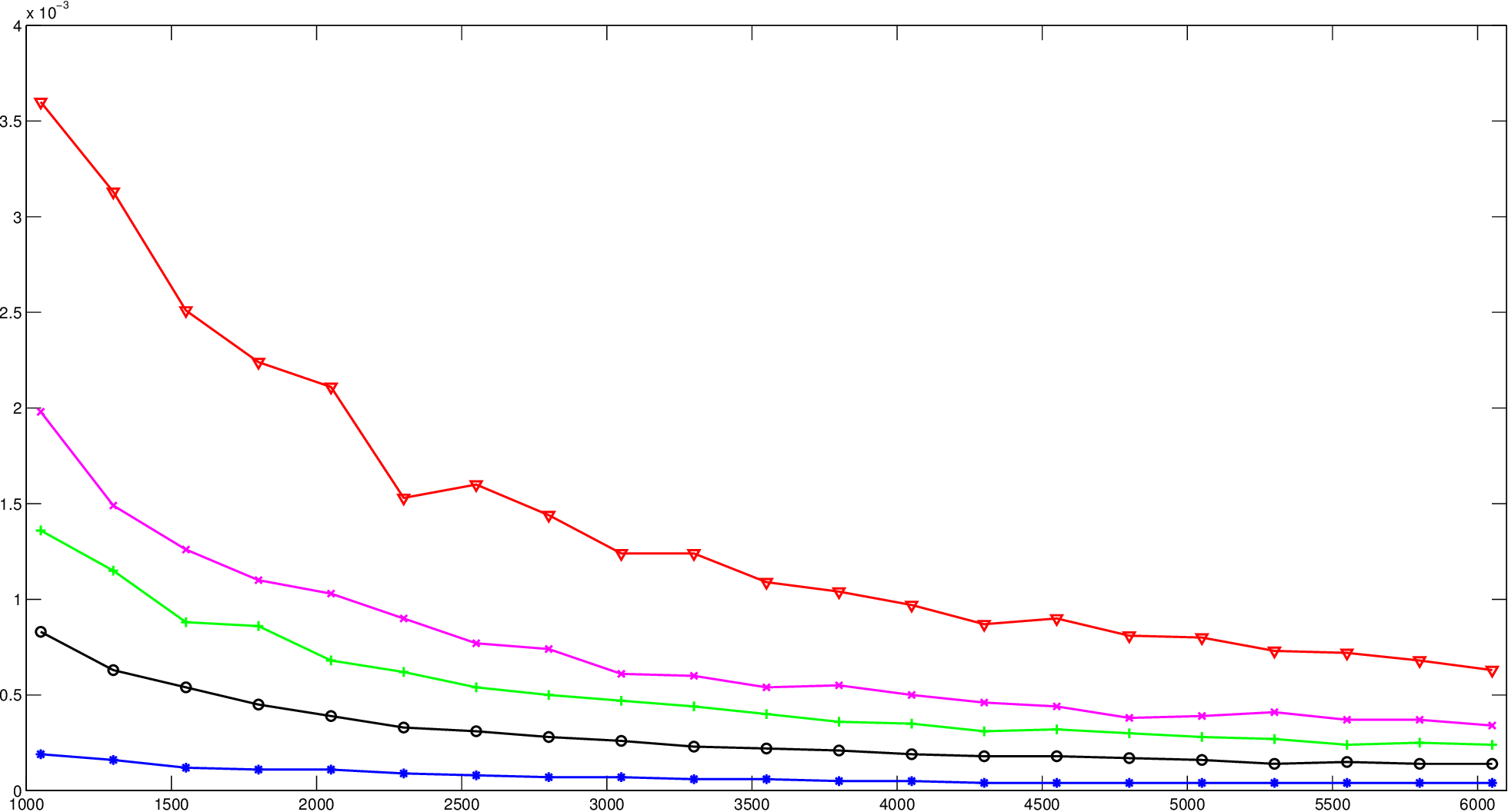}
  \label{A2:fig:Imagen4}
  
  \vspace{-0.1cm}
  \caption[\hspace{0.7cm} Empirical mean quadratic errors on the estimation of $ \theta$ of an O.U. process.]{{\small ${\rm EMSE} (N,T,\theta )$ based on $N=1000$ generations of O.U. process, for different sample sizes and values $\theta
  =0.1$ (blue star line), $\theta
  =0.4$ (black circles line), $\theta
  =0.7$ (green plus line), $\theta = 1$ (magenta cross line) and $\theta= 2$ (red triangle line).}}
   \label{A2:fig:3cc}
\end{figure}

Regarding asymptotic efficiency stated in \textcolor{Crimson}{Theorem} \ref{A2:th3prem},
from $N=1000$ simulations of the O.U. process over the interval
$[0,T],$ for $\left\lbrace T_l = 50 + 250(l-1), \ l=1,\dots,25\right\rbrace,$ the corresponding
empirical mean square errors $${\rm EMSE} (N,T,\theta ) = \frac{1}{N}
\displaystyle \sum_{s=1}^{N} \left(\theta -
\widehat{\theta}_{T}(\omega_{s}) \right)^2, \quad N=1000, \quad \theta =  \left[0.1, 0.4, 0.7, 1\right],$$  are displayed in Figure
\ref{A2:fig:3cc}. Here, $\widehat{\theta}_{T}(\omega_{s})$, with
$\omega_{s}\in \Omega,$ $s=1,\dots,N,$ represent the respective
approximated values (\ref{A2:173}) of the MLE of $\theta ,$ computed
from $\xi_{t_{i},s},$ $s=1,\dots,N,$ $t_{i}\in [0,T],$
$i=1,\dots,n.$ It can be observed, from the results displayed in
Figure \ref{A2:fig:3cc}, that \textcolor{Crimson}{Theorem} \ref{A2:th3prem}  holds for $T$
sufficiently large.

\textcolor{Crimson}{\subsection{Consistency of  $\rho_{\widehat{\theta}_{T}}=\rho_{\widehat{\theta}_{n}}$ in $\mathcal{L}(H)$ and $\mathcal{L}(B)$}
\label{A2:sec:92}}

The strong--consistency of $\rho_{\widehat{\theta}_n}$
in $\mathcal{L}(H)$ is derived in  \textcolor{Crimson}{Proposition} \ref{A2:prf} from    the
following almost surely upper bound

\begin{equation}
\Vert \rho_{\theta} - \rho_{\widehat{\theta}_n} \Vert_{\mathcal{L}
\left(H \right)} \leq   \vert \theta
-\widehat{\theta}_{n}\vert h \sqrt{\frac{h}{3} + 1} \quad a.s. \label{A2:177}
\end{equation}

\begin{figure}[H]
\centering
    \includegraphics[width=0.8\textwidth]{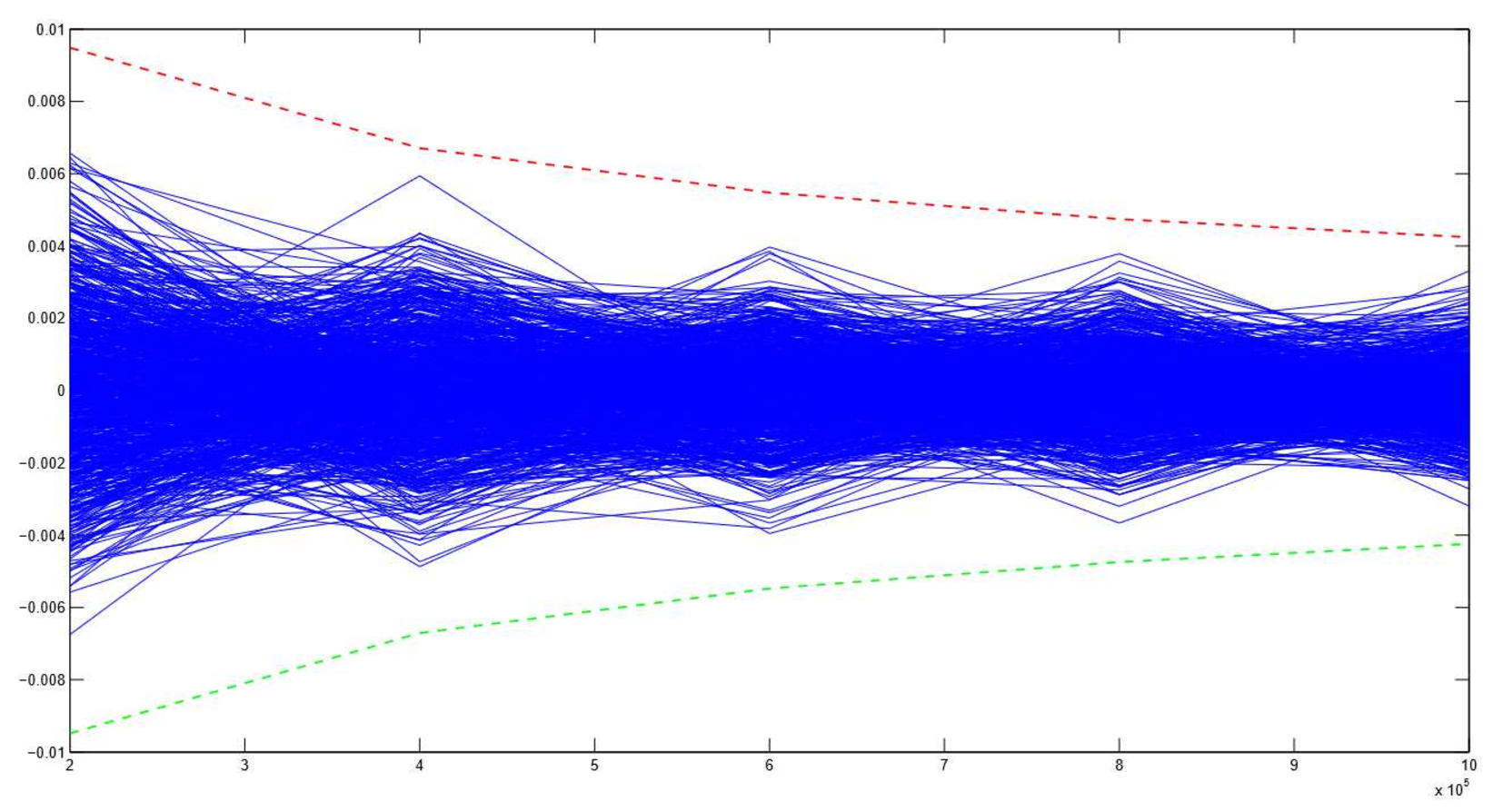}
  \label{A2:fig:Imagen4}

\vspace{-0.1cm}
\caption[\hspace{0.7cm} Empirical absolute errors on the estimation of $\theta$ of an O.U. process when large sample sizes are considered.]{{\small The values of  $\widehat{\theta}_{T} - \theta $ are
represented,  corresponding to $N = 1000$  simulations of  an O.U.
process over the interval $[0,T],$ with $ \left\lbrace T_l= 200000+(l-1)200000, \ l=1,\dots, 5 \right\rbrace,$ considering  the cases $\theta = 0.4$ (at the top), and $\theta = 1$ (at the bottom).
 The upper red dotted line is
$+3\sqrt{\frac{2 \theta}{T}}$  and the lower green dotted line is
$-3\sqrt{\frac{2 \theta}{T}}$.}} \label{A2:fig:2aa}
\end{figure}

 Here,  from $N=1000$ simulations of the O.U. process on
the interval $[0,T],$ with sample sizes \linebreak $T=nh=n= \left\lbrace 200000+(l-1)200000, \ l=1,\dots, 5 \right\rbrace,$  the corresponding values of $\widehat{\theta}_{T} -
\theta =\widehat{\theta}_{n}- \theta$ are computed,  considering the
cases $\theta = \left[0.4,  0.7, 1 \right].$  Table \ref{A2:tab:Table2baba} shows the empirical probability of
$\widehat{\theta}_{T} - \theta $ to lie within the band $\pm
3\sqrt{\frac{2 \theta}{T}},$ for each one of  sample sizes  and cases $\theta = \left[ 0.4, 0.7,
1 \right]$ regarded. It can be observed that for the sample sizes
studied, in the case of $\theta = 1,$ the empirical probabilities
are equal to one. Thus, the almost surely convergence to zero of the
upper bound (\ref{A2:177}) holds, with an approximated  convergence rate of
$\sqrt{T}=\sqrt{n}.$   Note that, for the other two cases, $\theta =
0.4$ and $\theta = 0.7,$ the empirical probabilities are also very
close to one (see also  Table \ref{A2:tab:Table1} for smaller sample
sizes, where
 we can also observe the empirical probabilities very close to one for
 the same band). In particular, Figure \ref{A2:fig:2aa}  displays the cases $\theta = 0.4$ (at the top) and $\theta =
 1$ (at the bottom).

\begin{table}[H]
\caption[\hspace{0.7cm} Empirical probability for the errors on the estimation of $\theta$ of an O.U. process for large sample sizes.]{{\small Empirical probability of $\widehat{\theta}_{T} - \theta$ to
be  within the band  $\pm 3\sigma = \pm 3\sqrt{\frac{2 \theta}{T}},$
from $N=1000$ simulations of an O.U. process over the interval $[0,T],$
with  $ \left\lbrace T_l= 200000+(l-1)200000, \ l=1,\dots, 5 \right\rbrace,$ considering the
cases $\theta = \left[ 0.4, 0.7, 1 \right]$.}}
\centering
\begin{small}
\begin{tabular}{|c||c|c|c|}
  \hline
   & \multicolumn{3}{c|}{Parameter $\theta$} \\
   \hline 
   $T $ & $0.4$  & $0.7$  & $1$ \\
  \hline \hline
   $200000$ & $1$ & $1$  & $1$  \\
  \hline
   $400000$ & $1$ & $1$ & $1$ \\
  \hline
     $600000$ & $0.999$ & $1$ & $1$ \\
  \hline
     $800000$ & $0.999$ & $0.999$ & $1$  \\
  \hline
     $1000000$ & $0.998$ & $1$ & $1$ \\
  \hline
\end{tabular}
\end{small}
  \label{A2:tab:Table2baba}
\end{table}

\begin{figure}[H]
   \centering
    \includegraphics[width=0.8\textwidth]{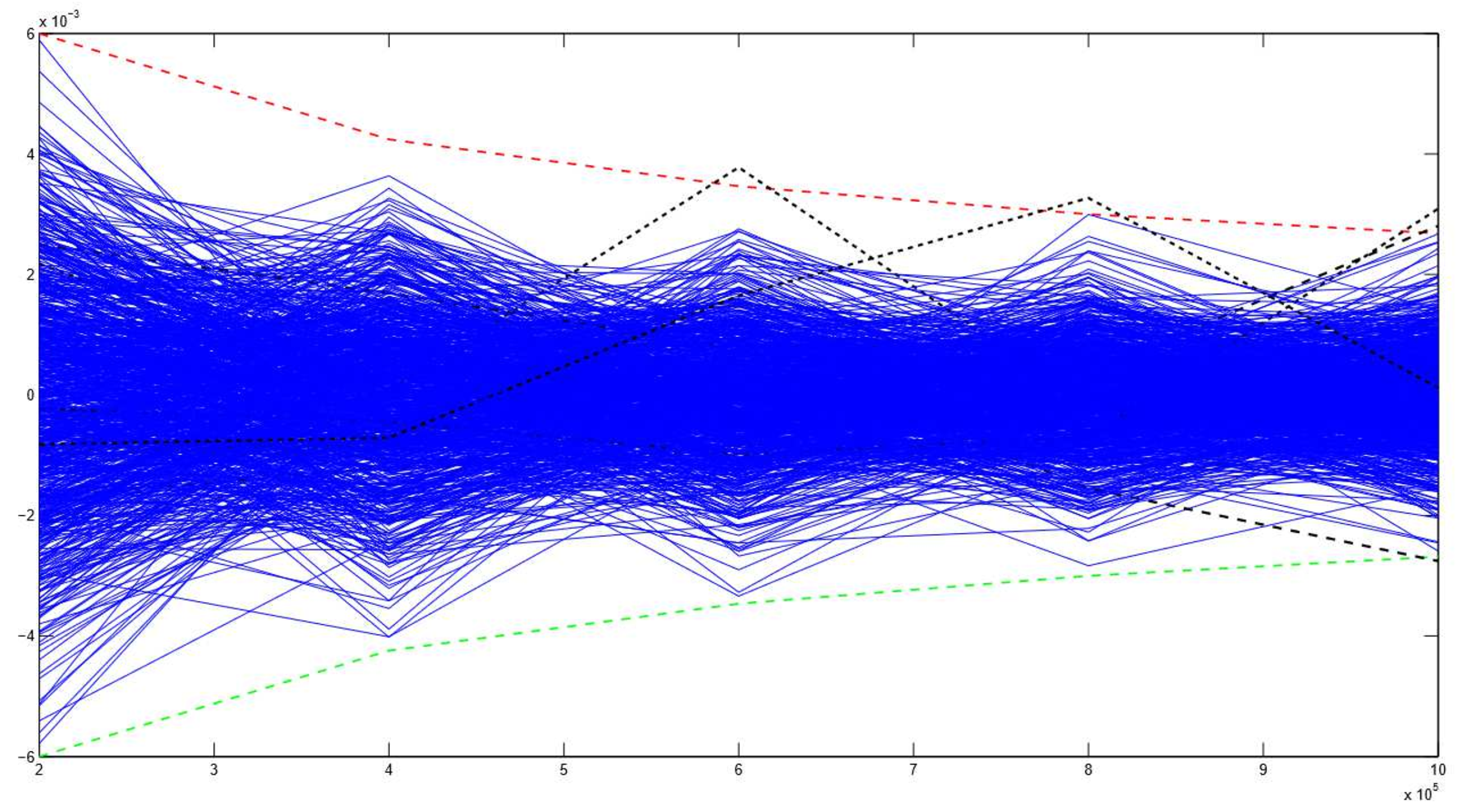}
  \label{A2:fig:Imagen3}
\end{figure}

It can be observed  from Table \ref{A2:tab:Table2baba} that  a better
performance is obtained for the largest values of  $\theta,$ which
corresponds to the weakest dependent case. Furthermore, from the
upper bound  in (\ref{A2:eqscbanach}), the strong consistency of
$\rho_{\widehat{\theta }_{n}}$ in $\mathcal{L}(B),$ with, as before,
$B=\mathcal{C}([0,h]),$ is also illustrated from the results
displayed in Table \ref{A2:tab:Table2baba} and Figure \ref{A2:fig:2aa}.

\textcolor{Crimson}{\subsection{Consistency of  the ARH(1) and  ARB(1) plug--in predictors for the O.U. process}
\label{A2:sec:93}}

\begin{figure}[H]
    \includegraphics[width=0.95\textwidth]{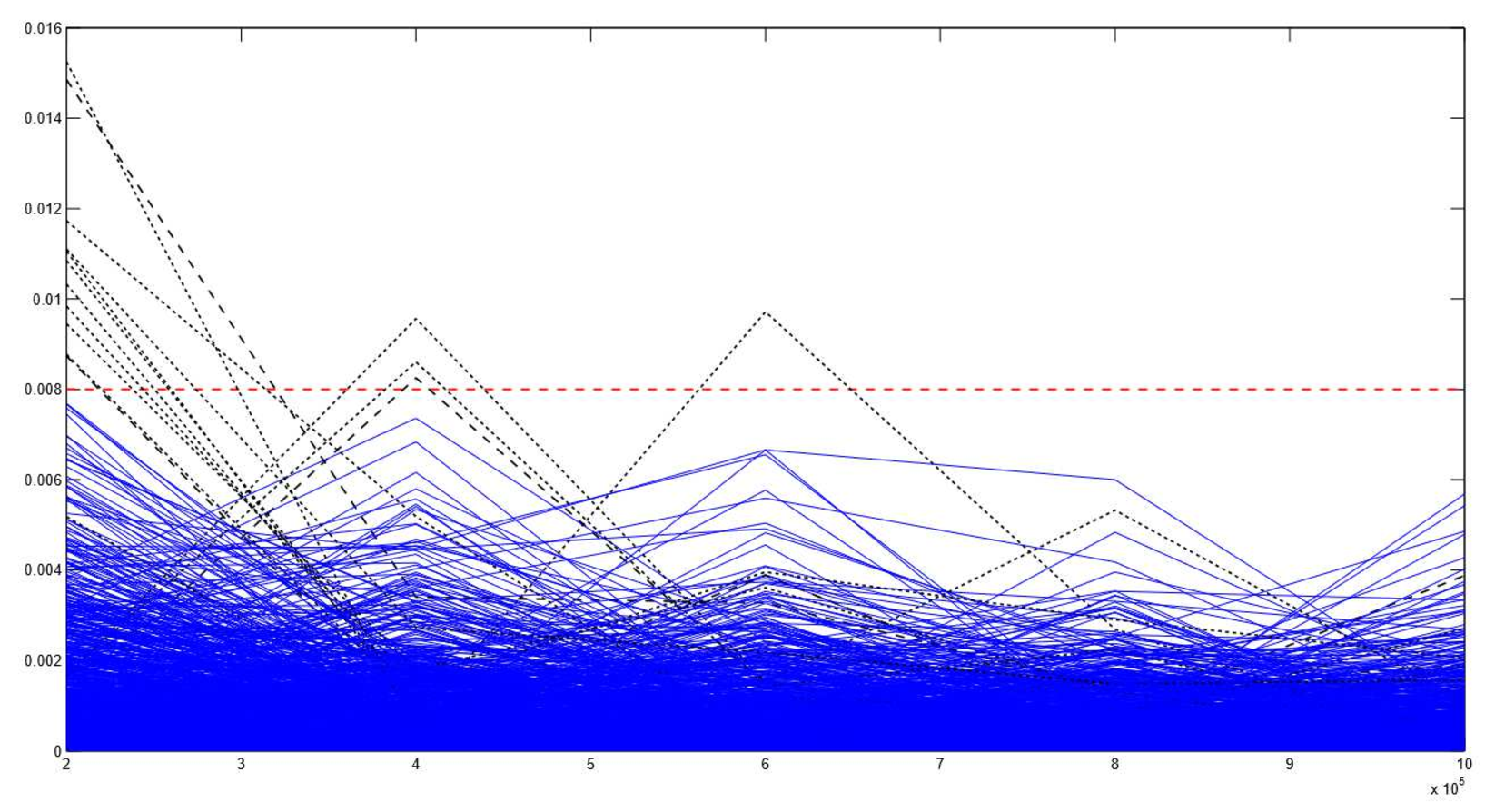}
  \label{A2:fig:Imagen3b}
  
\vspace{-0.1cm}
    \includegraphics[width=0.95\textwidth]{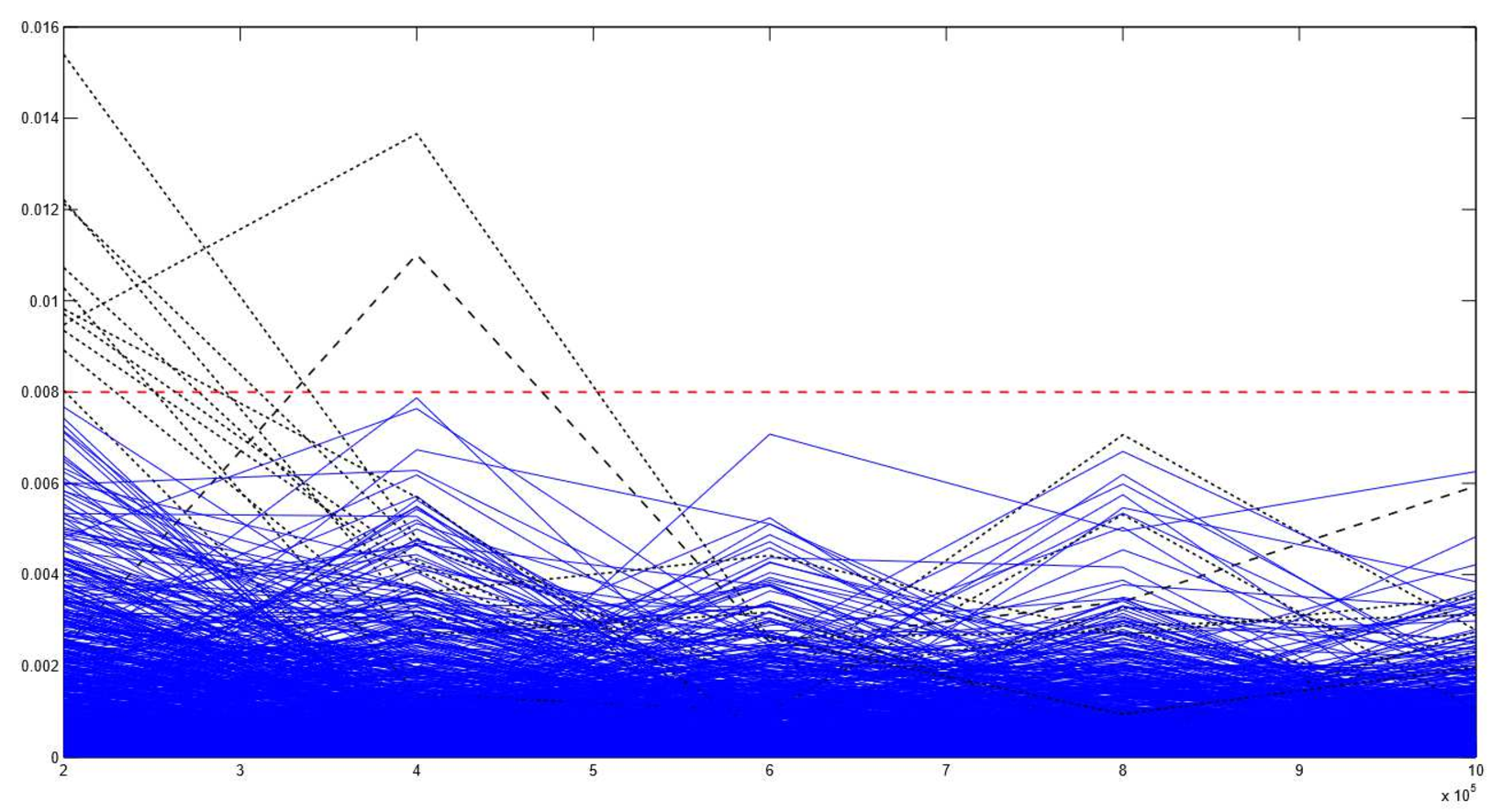}
  \label{A2:fig:Imagen4b}
  
\vspace{-0.1cm}
\caption[\hspace{0.7cm} Consistency of the ARH(1) and ARB(1) plug--in predictors for the O.U. process for different parameters $\theta$.]{{\small The values  of $\vert X_{n-1}\left( h \right) \vert \vert
\theta -\widehat{\theta}_n \vert h \sqrt{\frac{h}{3} + 1}$ (first and second figure) and
$\vert X_{n-1}\left( h \right) \vert \vert \theta
-\widehat{\theta}_n \vert h$ (third and fourth figure) are represented,  based on $N
= 1000$ generations  of O.U. process over the interval $[0,T],$ for
$ \left\lbrace T_l = 200000+(l-1)200000, \ l=1,\dots, 5 \right\rbrace,$ against $\epsilon =
0.008$ (red dotted line), considering $\theta = 0.4$ (first and third figure) and $\theta = 1$ (second and fourth figure).}}
\label{A2:fig:3}
\end{figure}

Let us now consider the derived upper bounds in (\ref{A2:53}) and
(\ref{A2:61})  in \textcolor{Crimson}{Corollaries} \ref{A2:c1}--\ref{A2:c2}, for the ARH(1) and
ARB(1) predictors, respectively. From the generation of  $N=1000$
discrete realizations of  an O.U. process over the interval
$[0,T],$ for $ \left\lbrace T_l=200000+(l-1)200000, \ l=1,\dots, 5 \right\rbrace,$ the upper
bounds (\ref{A2:53}) and (\ref{A2:61}) are evaluated, for the cases
$\theta = \left[ 0.4, 0.7, 1 \right].$ The following empirical probabilities for
$\epsilon=0.008,$ are  reflected in Table
\ref{A2:tab:Table3}

\begin{eqnarray}
\widehat{\mathcal{P}}^{H}(N,T,\theta)&=&1-\widehat{\mathcal{P}}\left(\vert X_{n-1}\left( h \right) \vert \vert
\theta -\widehat{\theta}_n \vert   h \sqrt{\frac{h}{3} +
1}>\epsilon\right),\label{A2:ep1H}\\ 
\widehat{\mathcal{P}}^{B}(N,T,\theta)&=&1-\widehat{\mathcal{P}}\left(\vert
X_{n-1}\left( h \right) \vert \vert \theta -\widehat{\theta}_n \vert
h >\epsilon\right), \label{A2:ep2B} 
\end{eqnarray}
\noindent with $N=1000$, $ \left\lbrace T_l=200000+(l-1)200000,\ l=1,\dots, 5 \right\rbrace$ and  $\theta = \left[ 0.4, 0.7, 1 \right]$,  for the  Hilbert--valued and
Banach--valued (see (\ref{A2:53}) and (\ref{A2:61})) frameworks (see also Figure
\ref{A2:fig:3}). It can be observed that the empirical probabilities
are equal to one in both frameworks for the largest sample sizes, in
any of the cases considered.

\begin{table}[H]
\caption[\hspace{0.7cm} Consistency of the ARH(1) and ARB(1) plug--in predictors for the O.U. process.]{\small{Empirical probabilities (\ref{A2:ep1H})--(\ref{A2:ep2B}),
based on $N=1000$ simulations of the O.U. process over the interval
$[0,T],$ for $ \left\lbrace T_l= 200000+(l-1)200000, \ l=1,\dots, 5 \right\rbrace,$ considering
the cases $\theta = \left[0.4, 0.7, 1\right],$  and $\epsilon=0.008$.}}
\centering
\begin{small}
\begin{tabular}{|c||ccc||ccc|}
  \hline
     &  \multicolumn{6}{c|}{Parameter $\theta$} \\
     \hline
& & Hilbert-valued case  & & & Banach-valued case & \\
\hline
   $T $ & $0.4$  & $0.7$  & $1$ & $0.4$  & $0.7$  & $1$\\
  \hline \hline
   $200000$ & $0.980$ & $0.980$ & $0.980$ & $0.987$ & $0.991$  & $0.987$ \\
  \hline
   $400000$ & $0.995$ & $0.995$ & $0.995$ & $0.997$ & $0.998$ & $0.9977$ \\
  \hline
     $600000$ & $0.999$ & $0.998$ & $0.999$  & $0.999$ & $0.999$ & $1$\\
  \hline
     $800000$ & $1$ & $0.999$ & $0.999$  & $1$ & $1$ & $1$\\
  \hline
     $1000000$ & $1$ & $1$ & $1$ & $1$ & $1$ & $1$ \\
  \hline
\end{tabular}
\end{small}
  \label{A2:tab:Table3}
\end{table}

\begin{figure}[H]
\centering
    \includegraphics[width=0.95\textwidth]{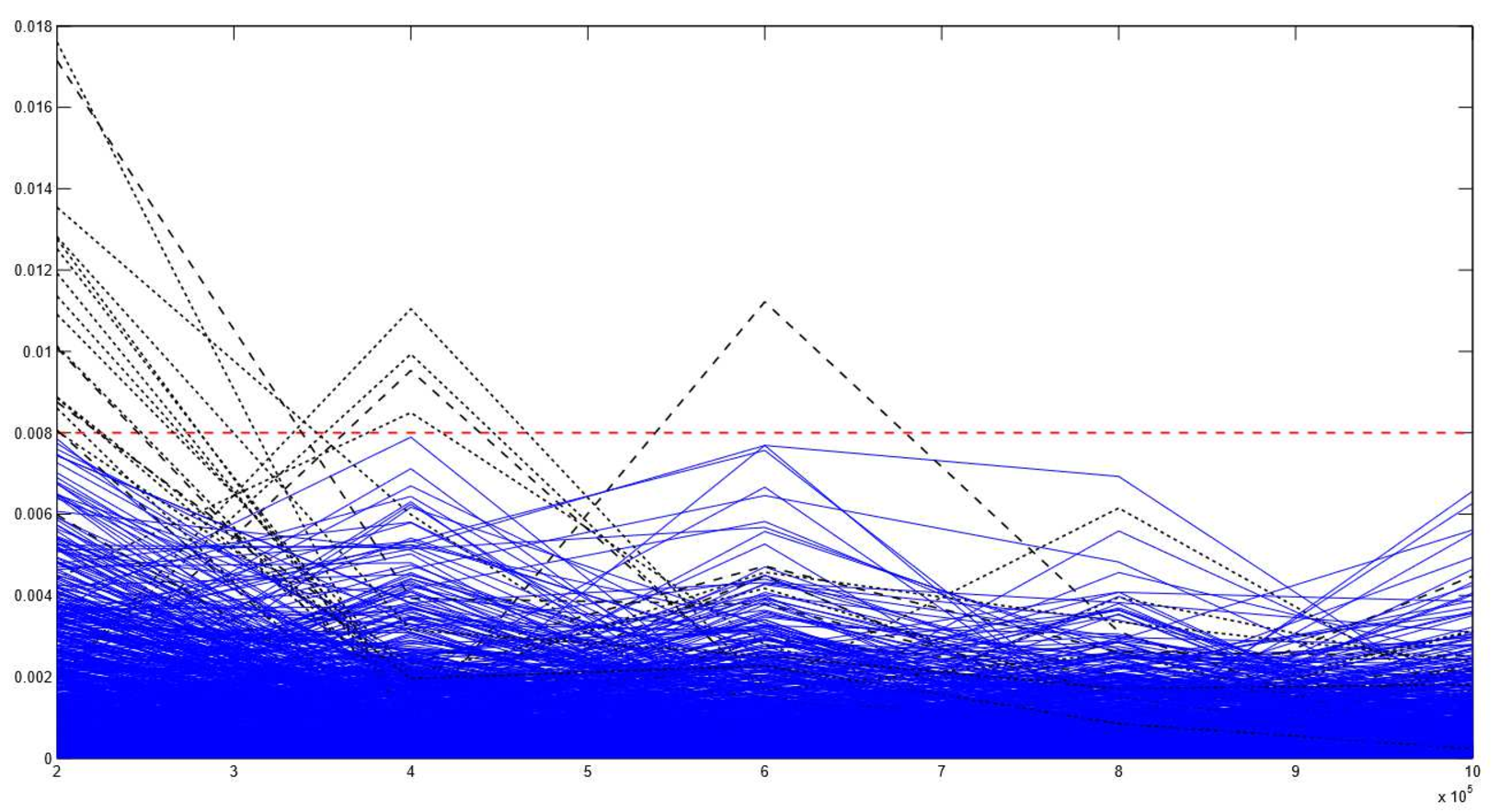}
  \label{A2:fig:Imagen3}

\vspace{-0.1cm}
    \includegraphics[width=0.95\textwidth]{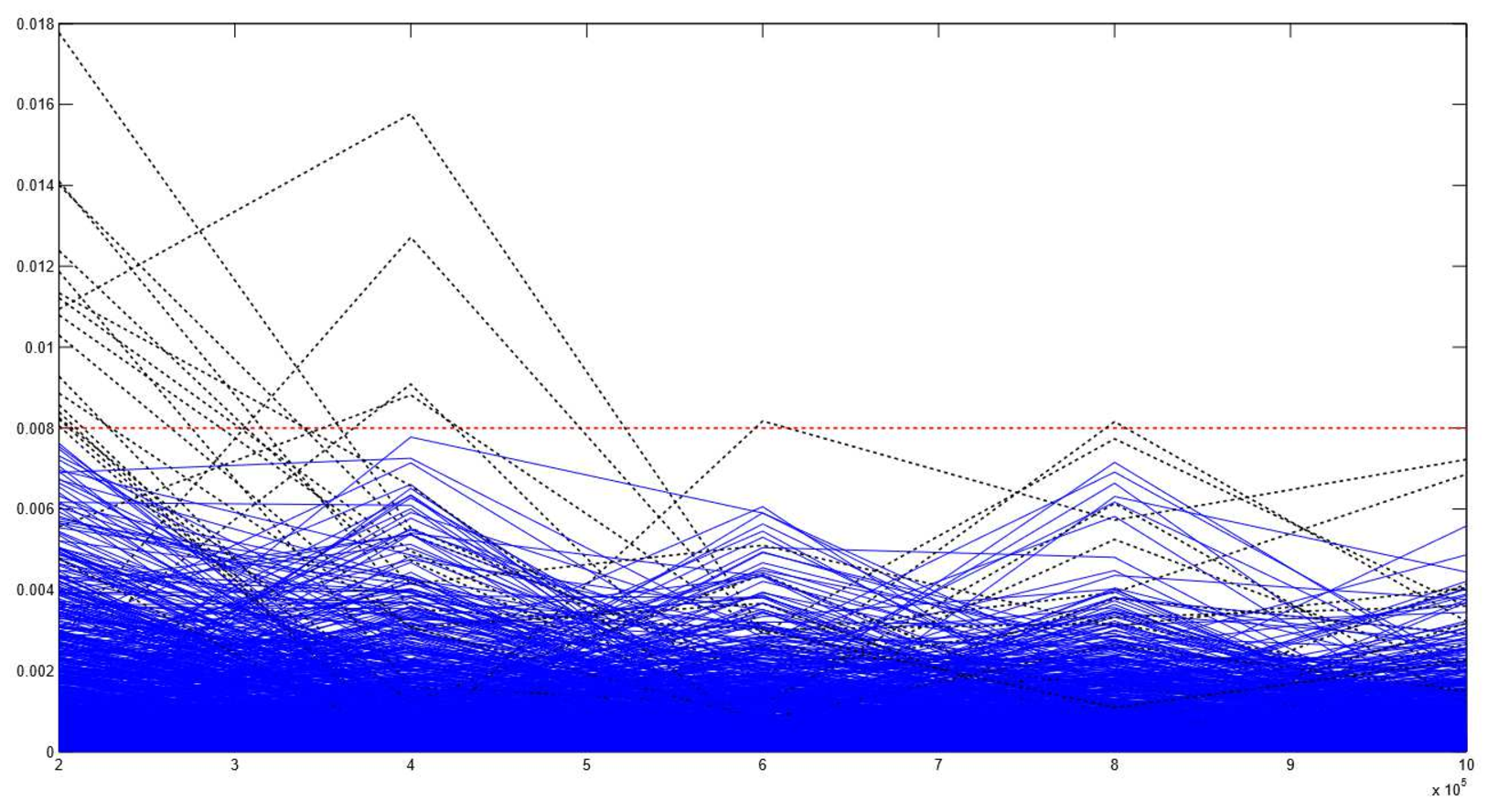}
  \label{A2:fig:Imagen4}
\end{figure}

 The
strong--consistency of the MLE of $\theta $ and of the
autocorrelation operator of the O.U. process, in Banach and Hilbert
spaces, has been first illustrated. The almost surely rate of
convergence to zero is shown as well. The numerical results on  the
consistency of the associated ARH(1) and ARB(1) plug--in predictors
then follow, from the computation of the corresponding empirical
probabilities for the derived upper bounds. Note that the numerical
results displayed in \textcolor{Crimson}{Appendix} \ref{A2:sec:9} are obtained under
generation of sample sizes ranging from $12000$ up to a million of
time instants, considering $1000$ repetitions for each one of such
sample sizes. In all these simulations performed, the discretization
step size considered has been $\Delta t=0.02.$

\textcolor{Crimson}{\section{Final comments}
\label{A2:sec:7}}

The problem of functional prediction of the O.U.
process could be of interest in several applied fields. For example,
in finance, in the context of the Vasicek's model (see \cite{Vasicek77})
the results derived allow to predict the curve representing the
interest rate over a temporal interval, in a consistent way. Note
that, in this context,  the MLE computed for parameter
$\theta $ provides a consistent approximation of the speed
reversion, which definitely determines the proposed  functional
predictor of the interest rate.

Summarizing, this paper addresses the problem of functional
prediction of the O.U. process from ARH(1) and ARB(1) perspectives.
Specifically, considering the O.U. process as an ARH(1) and an
ARB(1) process, new results on strong consistency (almost surely
convergence to the true parameter value), in the spaces
$\mathcal{L}(H)$ and $\mathcal{L}(B)$  of the MLE  of its
autocorrelation operator are derived. Consistency results
(convergence in probability to the true value) of the associated
plug--in predictors are obtained as well. The numerical results
shown, in addition,  the normality and the asymptotic efficiency of
the MLE of the scale parameter $\theta $ of the covariance function
of the O.U. process.

\textcolor{Crimson}{\section{Supplementary Material}
\label{A2:Supp}}

The definition and properties of an O.U. process are given here, as well as the proof of \textcolor{Crimson}{Lemma} \ref{A2:lem1}.

\textcolor{Crimson}{\subsection{Ornstein--Uhlenbeck process} \label{A2:Supp:sec:21}}

Let $\xi \left(\omega \right)= \left\lbrace \xi_t \left(\omega
\right), \ t \in \mathbb{R} \right\rbrace,$ $\omega \in \Omega ,$ be a
real--valued sample--path continuous stochastic process defined on the
basic probability space $\left(\Omega,\mathcal{A},\mathcal{P}\right),$
with index set the real line $\mathbb{R}.$  As demonstrated in
\cite{Doob42}, process $\xi $ is an O.U. process if it provides  the
Gaussian solution to the following stochastic linear Langevin
differential equation:
\begin{equation}
d \xi_t = \theta \left( \mu - \xi_t \right)dt + \sigma
dW_t,\quad\theta,~\sigma > 0,\quad t \in \mathbb{R}, \label{A2:Supp:4}
\end{equation}
\noindent where $W = \left \lbrace W_t, \ t \in \mathbb{R} \right\rbrace$ is a
standard bilateral Wiener process; i.e., $$W_t = W_{t}^{(1)}
\boldsymbol{1}_{\mathbb{R}^+} \left( t \right) + W_{-t}^{(2)}\boldsymbol{1}_{\mathbb{R}^-}\left( t
\right),$$ with $W_{t}^{(1)}$ and $W_{-t}^{(2)}$ being  independent
standard Wiener processes, and $\boldsymbol{1}_{\mathbb{R}^+}$ and   $\boldsymbol{1}_{\mathbb{R}^-}$
respectively denoting the indicator functions over the positive and
negative real line. Applying, in equation  (\ref{A2:Supp:4}), the method of
separation of variables, considering $f\left(\xi_t,t\right) = \xi_t
e^{\theta t},$ we obtain
\begin{eqnarray}
\xi_t &=& \mu + \displaystyle \int_{-\infty}^{t}\sigma e^{-\theta
\left(t - s \right)} dW_s, \quad \theta, \quad \sigma > 0, \quad t \in \mathbb{R},
\label{A2:Supp:5}
\end{eqnarray}
\noindent where the integral is understood in the It\^{o}
sense (see \cite{AshGardner76,Sobczyk91} for more details).
Particularizing to $\xi = \left \lbrace \xi_t, \ t \in
\mathbb{R}^+ \right\rbrace$, the O.U. process is transformed into
\begin{equation}
\xi_t = \xi_0 e^{-\theta t}  + \mu\left( 1 - e^{-\theta t} \right) +
\displaystyle \int_{0}^{t}\sigma e^{-\theta \left(t - s \right)}
dW_s,\quad \theta,\sigma > 0, \quad t \in \mathbb{R}^+. \label{A2:Supp:6}
\end{equation}

It is well--known that the solution  $\xi = \left \lbrace \xi_t, \ t \in \mathbb{R} \right\rbrace$ to the stochastic differential equation
\begin{equation}
d\xi_t = \mu \left(\xi_t,t\right)dt +
\sqrt{D\left(\xi_t,t\right)}dW_t, \quad t \in \mathbb{R}, \nonumber % \label{A2:Supp:7}
\end{equation}
\noindent has  marginal probability density function $f\left(x,t \right),$
satisfying the following Fokker--Planck's scalar equation (see, for
example, \cite{Kadanoff00}):
\begin{equation}
\frac{\partial }{\partial t} f\left(x,t\right) =
\frac{-\partial}{\partial x}
\left[\mu\left(x,t\right)f\left(x,t\right) \right] + \frac{1}{2}
\frac{\partial^2}{\partial x^2} \left[D\left(x,t\right)
f\left(x,t\right)\right], \quad t \in \mathbb{R}. \nonumber %\label{A2:Supp:8}
\end{equation}

In the case of O.U. process, the stationary solution
($\frac{\partial}{\partial t} f\left(x,t\right) = 0$), under
$f\left(x,x_0\right) =  \delta\left(x - x_0 \right)$, adopts the
form
\begin{equation}
f\left(x,t\right) = \sqrt{\frac{\theta}{\pi \sigma^2}}
e^{\frac{-\theta \left(x - \mu \right)^2}{\sigma^2}}, \quad
\theta,\sigma > 0, \quad t \in \mathbb{R}, \nonumber %\label{A2:Supp:9}
\end{equation}
\noindent which corresponds to the probability density function of a
Gaussian distribution with mean $\mu $ and variance
$\frac{\sigma^2}{2\theta},$  i.e., which corresponds to the
probability density function of a random variable $X$ such that
$$X\sim\mathcal{N}\left(\mu,\frac{\sigma^2}{2\theta}
\right).$$ From (\ref{A2:Supp:5}), the mean  and covariance functions of O.U.
process (see, for instance, \cite{Doob42,UhlenbeckOrnstein30}) can be
computed as follows:
\begin{eqnarray}
\mu_{\xi}(t)&=& {\rm E} \left\lbrace \xi_t\right\rbrace=  \mu + \sigma {\rm E} \left\lbrace \displaystyle \int_{-\infty}^{t} e^{-\theta \left(t - s \right)} dW_{s}\right\rbrace =\mu,\quad t \in \mathbb{R}, \nonumber \\
C_{\xi}(t,s)&=& {\rm Cov} \left(\xi_s,\xi_t\right)  =  {\rm E}\left\lbrace\left(\xi_s - \mu \right)\left(\xi_t - \mu\right) \right\rbrace = \sigma^2 e^{-\theta \left(t + s \right)} {\rm E} \left\lbrace \displaystyle \int_{-\infty}^{t} e^{\theta u} dW_u\displaystyle \int_{-\infty}^{s} e^{\theta v} dW_v\right\rbrace  \nonumber \\
&=& \sigma^2 e^{-\theta \left(t + s \right)} \int_{-\infty}^{\infty} e^{2\theta u}\boldsymbol{1}_{\left[-\infty,t\right]}\left(u\right) \boldsymbol{1}_{\left[-\infty,s\right]}\left(u\right) du = \sigma^2 e^{-\theta \left(t + s \right)} \int_{-\infty}^{\min\left\lbrace s,t\right\rbrace} e^{2 \theta u}du  \nonumber \\
&=&\frac{\sigma^2}{2\theta} e^{-\theta \left(t + s \right)}
e^{2\theta\min\left\lbrace s,t \right\rbrace}  =
\frac{\sigma^2}{2\theta} e^{-\theta \vert t - s \vert},\quad t,s \in
\mathbb{R},  \label{A2:Supp:10}
\end{eqnarray}
\noindent where ${\rm Cov} \left(X,Y \right)$ denotes the
covariance between random variables $X$ and $Y$. Additionally, from
(\ref{A2:Supp:6}), we obtain the following identities:
\begin{eqnarray}
{\rm E} \left\lbrace\xi_t\right\rbrace &=& \mu e^{-\theta t}  + \mu\left( 1 - e^{-\theta t} \right) = \mu, \quad {\rm E}\left\lbrace \xi_t \vert \xi_0 = c\right\rbrace = \mu  + e^{-\theta t} \left(c - \mu \right),\quad t \in \mathbb{R}^+,\nonumber \\
{\rm Cov} \left(\xi_s,\xi_t \vert  \xi_0 = c\right) &=&
\frac{\sigma^2}{2 \theta}  e^{-\theta \vert t-s \vert} +\left(c^2 -
2c \mu + \mu^2 \right)e^{-\theta \left(s + t \right)},\quad t,s \in
\mathbb{R}^+, \nonumber %\label{A2:Supp:11}
\end{eqnarray}
\noindent where $c$ is a constant. In the subsequent development, we will
consider $\mu = 0$ and $\sigma = 1$.

\textcolor{Crimson}{\subsection{Maximum likelihood estimation of the covariance scale parameter $\theta$}
\label{A2:Supp:sec:23}}

The MLE of  $\theta$ in (\ref{A2:Supp:10}) is given by  (see
\cite{GraczykJakubowski06}; \cite[p. 63]{Kutoyants04}; \cite[p. 265]{LiptserShiraev01})
%\begin{equation}
%L\left(\widehat{\theta}; \xi_{(T)}\right) = e^{-\displaystyle \int_{0}^{T} S\left(\xi_t,t \right) d \xi_t + \frac{1}{2} \displaystyle \int_{0}^{T}S^{2}\left(\xi_t,t \right) dt},~T > 0, \label{23}
%\end{equation}
%where $S\left(\xi_t,t \right) = -\widehat{\theta} \xi_t$. From (\ref{5}) and maximizing the \textbf{log-likelihood function}, we get
\begin{equation}
\widehat{\theta}_{T} = \frac{-\displaystyle \int_{0}^{T} \xi_t d
\xi_t}{\displaystyle \int_{0}^{T}\xi_{t}^{2} dt} = \frac{\theta
\displaystyle \int_{0}^{T} \xi_{t}^{2} dt - \displaystyle
\int_{0}^{T}\xi_t d W_t}{\displaystyle \int_{0}^{T}\xi_{t}^{2} dt} =
\theta - \frac{ \displaystyle \int_{0}^{T}\xi_t d W_t}{\displaystyle
\int_{0}^{T}\xi_{t}^{2} dt},\quad \theta,T > 0. \label{A2:Supp:14}
\end{equation}

Thus, equation (\ref{A2:Supp:14}) becomes
\begin{eqnarray}
\widehat{\theta}_{T} &=& \frac{1 + \frac{\xi_{0}^{2}}{T} -
\frac{\xi_{T}^{2}}{T}}{\frac{2}{T} \displaystyle \int_{0}^{T}\xi_{t}^{2} dt}, \quad
T > 0. \label{A2:Supp:17}
\end{eqnarray}

We will assume that $T$ is large enough such that
$\widehat{\theta}_T
> 0$ almost surely. It is well--known that the MLE $\widehat{\theta}_T$ of $\theta$   is
strongly consistent (see details in \cite[Proposition 2.2]{KleptsynaBreton02}; \cite[p. 63 and p. 117]{Kutoyants04}).

\begin{theorem}
\label{A2:Supp:th2prem}\textit{ The following limit in distribution sense holds for
the MLE $\widehat{\theta}_T$ of $\theta,$ given in equation
(\ref{A2:Supp:17}):
\begin{equation}
\lim_{T\rightarrow \infty}\sqrt{T} \left(\widehat{\theta}_T-\theta
\right) = \lim_{T\rightarrow \infty }\frac{- \sqrt{T}
\displaystyle \int_{0}^{T}\xi_t d W_t}{ \displaystyle \int_{0}^{T}\xi_{t}^{2}
dt}=Z,\quad \mbox{with}\quad  Z\sim \mathcal{N} \left(0, 2 \theta
\right). \nonumber %\label{A2:Supp:19}
\end{equation}}
\end{theorem}

Results in \cite[Theorem 1.1 and Corollary 1.1]{Jiang12} lead to the
following almost surely identities (see also \cite[Theorem 2.10]{Bosq00};\cite[pp. 196--203]{LedouxTalagrand91}, in relation to the
law of the iterated logarithm)
\begin{eqnarray}
&&\displaystyle  \limsup_{T \to +\infty} \frac{  \widehat{\theta}_T
-\theta}{\sqrt{\frac{4 \theta}{T} \displaystyle \ln \left(
\displaystyle \ln \left(T \right)\right)}} = 1\quad
a.s., \nonumber \\ %\label{A2:Supp:19b}\\
&&\displaystyle  -\liminf_{T \to +\infty} \frac{ \widehat{\theta}_T
-\theta}{\sqrt{\frac{4 \theta}{T} \displaystyle \ln \left(
\displaystyle \ln \left(T \right)\right)}}=1 \quad
a.s.,
\nonumber \\ %\label{A2:Supp:20}\\
&&\vert \theta - \widehat{\theta}_T \vert = \mathcal{O}
\left(\sqrt{\frac{4 \theta \displaystyle \ln\left( \displaystyle
\ln \left(T \right)\right)}{T}}\right) \quad
a.s. \nonumber %\label{A2:Supp:21}
\end{eqnarray}

\textcolor{Crimson}{\subsection{Preliminary inequalities and results}
\label{A2:Supp:sec:4}}

In this section we recall some inequalities and well--known
convergence  results on random variables, as well as basic
deterministic inequalities,  that have been  applied in the
derivation of the main results displayed above.

\begin{lemma}
\label{A2:Supp:lempr1} \textit{Let $X$ be a zero--mean normal distributed random
variable, i.e.,  $X \sim \mathcal{N} \left(0, \sigma^2 \right),$
with $\sigma > 0$. Then,
\begin{equation}
\mathcal{P} \left( \vert X \vert \geq x \right) \leq
e^{-\frac{x^2}{2 \sigma^2}}, \quad  x \geq 0. \nonumber %\label{A2:Supp:30}
\end{equation}}
\end{lemma}

\begin{proof}
Let $X^{\prime }$ be such that  $X^{\prime } \sim \mathcal{N}
\left(0, 1 \right).$ Then,
\begin{equation}
\mathcal{P} \left( \vert X^{\prime } \vert \geq x \right) =
2F_{X^{\prime }} \left(-x \right) =   \sqrt{\frac{2}{\pi}}
\displaystyle \int_{x}^{\infty} e^{- \frac{t^2}{2}}dt,~\forall x
\geq 0. \label{A2:Supp:31}
\end{equation}

Let us  set
\begin{eqnarray}
g \left( x \right) &=& e^{-\frac{x^2}{2}} - \sqrt{\frac{2}{\pi}}
\displaystyle \int_{x}^{\infty} e^{- \frac{t^2}{2}}dt, \quad
g\left(0 \right) = 0,\quad
\lim_{x\rightarrow \infty}g\left(x \right)=0, \nonumber \\
g' \left( x \right) &=& - xe^{-\frac{x^2}{2}}
+\sqrt{\frac{2}{\pi}}e^{-\frac{x^2}{2}} = e^{-\frac{x^2}{2}} \left(
\sqrt{\frac{2}{\pi}} - x \right). \nonumber
\\ \label{A2:Supp:32}
\end{eqnarray}
Function $g$ is monotone increasing   over $\left(0,
\sqrt{\frac{2}{\pi}} \right),$ and $g$ is monotone decreasing   over
$\left( \sqrt{\frac{2}{\pi}}, \infty \right).$

 From equations (\ref{A2:Supp:31})--(\ref{A2:Supp:32}),
$$\mathcal{P} \left( \vert X^{\prime } \vert \geq x \right) \leq
e^{-\frac{x^2}{2}}, \quad x \geq 0.$$ Now, consider
$X^{\prime} =\frac{X}{\sigma},$ with $X\sim \mathcal{N} \left(0,
\sigma^2 \right),$ then,
\begin{equation}
\mathcal{P} \left( \vert X \vert \geq x \right) \leq
e^{-\frac{x^2}{2 \sigma^2}}, \quad x \geq 0. \nonumber %\label{A2:Supp:33}
\end{equation}

\hfill \hfill \textcolor{Aquamarine}{$\blacksquare $}
\end{proof}

\textcolor{Crimson}{\subsubsection{Proof of Lemma 1}}

\begin{proof}

 Let us first consider the case $k=1,$  from
\begin{eqnarray} \rho_{\theta} \left(x\right)\left(t\right) &=&
e^{-\theta t} x\left(h \right), \quad \rho_{\theta} \left(X_{n-1}
\right)\left(t\right) = e^{-\theta t} \int_{-\infty}^{nh} e^{-\theta
\left(nh - s \right)} dW_s, \nonumber \\
 \varepsilon_{n} \left(t\right) &=&
\int_{nh}^{nh+t} e^{-\theta \left(nh + t - s \right)} dW_s, \nonumber %\label{A2:Supp:40}
\end{eqnarray}
\noindent and
\begin{eqnarray}\Vert \rho_{\theta}(x) \Vert_{H}^{2}  &=& \displaystyle \int_{0}^{h}
\left( \rho_{\theta} \left(x\right)\left(t\right) \right)^2 d
\left(\lambda + \delta_{(h)}\right)\left(t\right) = \displaystyle
\int_{0}^{h} \left( \rho_{\theta} \left(x\right)\left(t\right)
\right)^2 dt +  \left( \rho_{\theta} \left(x\right)\left(h\right)
\right)^2, \nonumber %\label{A2:Supp:42}
\end{eqnarray}

\noindent we have
\begin{eqnarray}
\Vert \rho_{\theta}  \Vert_{\mathcal{L}\left(H\right)} &=&
\displaystyle \sup_{x\in H} \left\lbrace \frac{\Vert
\rho_{\theta}\left(x \right) \Vert_H}{\Vert x \Vert_H} \right\rbrace
= \sup_{x\in H} \left\lbrace\sqrt{\frac{\left(\displaystyle
\int_{0}^{h}  e^{-2\theta t} dt + e^{-2\theta h}\right)
\left( x\left(h\right) \right)^2}{\displaystyle \int_{0}^{h} \left( x \left(t \right) \right)^2 dt
+ \left( x \left(h\right) \right)^2}}\right\rbrace. \label{A2:Supp:43}
 \label{A2:Supp:pr1}
\end{eqnarray}

 Furthermore, 
\begin{equation}\Vert \rho_{\theta}  \Vert_{\mathcal{L}\left(H\right)}=\sup_{x\in H} \left\lbrace \sqrt{\frac{\left(\displaystyle \int_{0}^{h}  e^{-2\theta t} dt +
e^{-2\theta h}\right) \left( x\left(h\right) \right)^2}{\displaystyle \int_{0}^{h}
\left( x \left(t \right)\right)^2 dt + \left( x \left(h\right) \right)^2}} \right\rbrace \leq \sqrt{\int_{0}^{h}
e^{-2\theta t} dt + e^{-2\theta h}}.\label{A2:Supp:eqow}
\end{equation} 

 Additionally, the function
$x_{0}:~[0,h]\longrightarrow \mathbb{R},$ given by
\begin{equation}x_{0}(t)=\chi_{\mathcal{M}}(t),\quad h\in \mathcal{M}\subset [0,h],\quad \int_{\mathcal{M}}dt=0, \label{A2:Supp:ex}
\end{equation} 
\noindent with $\boldsymbol{1}_{\mathcal{M}},$ denoting the indicator function of set $\mathcal{M},$ belongs to $H=L^2
\left(\left[0,h\right],\beta_{\left[0,h\right]},\lambda +
\delta_{(h)}\right),$ since $$x_{0}^{2}(h)=1, \quad \int_{0}^{h}x_{0}^{2}(t)dt=0 \quad \|x_{0}\|_{H}^{2}=\int_{0}^{h}x_{0}^{2}(s)ds+x_{0}^{2}(h)=1.$$

Thus, by definition of $\|\rho_{\theta
}\|_{\mathcal{L}(H)},$
\begin{eqnarray}\frac{\|\rho_{\theta }(x_{0})\|_{H}}{\|x_{0}\|_{H}}=\sqrt{\int_{0}^{h} e^{-2\theta t} dt + e^{-2\theta h}}&\leq & \|\rho_{\theta
}\|_{\mathcal{L}(H)}\label{A2:Supp:rw}
\end{eqnarray}

Equations (\ref{A2:Supp:43})--(\ref{A2:Supp:rw}) lead to
\begin{equation}\Vert \rho_{\theta}
\Vert_{\mathcal{L}\left(H\right)}= \sqrt{\int_{0}^{h} e^{-2\theta t}
dt + e^{-2\theta h}}=\sqrt{\frac{1 + e^{-2\theta h}\left(2\theta -
1\right)}{2\theta}}.\label{A2:Supp:fr}
\end{equation}

 We are now going to compute $\Vert
\rho_{\theta}^{k}\Vert_{\mathcal{L}\left(H\right)},$ for $k\geq 2.$
Since, for all $x\in H,$
\begin{equation}
\rho_{\theta}^{k}(x)(t)= e^{-\theta t}e^{-\theta (k-1)h}x(h), \nonumber %\label{A2:Supp:powerrho}
\end{equation}
\noindent we obtain
\begin{equation}
\Vert\rho_{\theta}^{k}\Vert_{\mathcal{L}\left(H\right)}=\sup_{x\in
H} \left\lbrace \sqrt{\frac{\left[e^{-2\theta (k-1)h} \displaystyle \int_{0}^{h}e^{-2\theta
t}dt+e^{-2\theta k
h}\right] \left(x(h) \right)^2}{\displaystyle \int_{0}^{h} \left( x(t) \right)^2 dt+ \left( x(h) \right)^2}} \right\rbrace. \nonumber %\label{A2:Supp:nr2}
\end{equation}

Considering function $x_{0}$ defined in equation (\ref{A2:Supp:ex}),
applying similar arguments to those given in the computation of
$\Vert \rho_{\theta}  \Vert_{\mathcal{L}\left(H\right)},$ we have
\begin{equation}
\Vert\rho_{\theta}^{k}\Vert_{\mathcal{L}\left(H\right)}=\sqrt{e^{-2\theta
(k-1)h}\frac{1+e^{-2\theta  h}\left(2\theta -1\right)}{2\theta}}=
e^{-\theta (k-1)h}\Vert \rho_{\theta}
\Vert_{\mathcal{L}\left(H\right)}. \nonumber %\label{A2:Supp:nr3}
\end{equation}

 Now, from equation (\ref{A2:Supp:fr}),

\begin{equation}
\Vert \rho_{\theta} \Vert_{\mathcal{L}\left(H\right)}<1
\Longleftrightarrow 1 - e^{-2\theta h} < 2\theta \left(1 - e^{-2\theta
h} \right) \Longleftrightarrow  \theta > \frac{1}{2}. \nonumber %\label{A2:Supp:45}
\end{equation}

Furthermore, for $\theta \in \left(0,1/2 \right]$,
\begin{equation}
\Vert \rho_{\theta} \Vert_{\mathcal{L}\left(H\right)} = \sqrt{\alpha
\left(\theta \right)} < \sqrt{1+h}, \nonumber
\end{equation}
\noindent since $\sqrt{\alpha \left(\theta \right)}$ is a
monotonically decreasing function on $\left(0,1/2 \right],$ with
$\alpha \left(\theta \right)=1$ if $\theta = \frac{1}{2}$ and \linebreak
$\alpha \left(\theta \right) \to 1+h$, when $\theta \to 0$. Hence,
if $\theta (k-1) \geq 1$,
\begin{equation}
\Vert \rho_{\theta}^{k} \Vert_{\mathcal{L}\left(H\right)}=
e^{-\theta(k-1)h}\sqrt{\alpha \left(\theta \right)} \leq e^{-h}
\sqrt{\alpha \left(\theta \right)}  < \frac{\sqrt{1+h}}{e^h} <
1,\quad h > 0, \nonumber
\end{equation}
\noindent which implies that $\Vert \rho_{\theta}^{k_0}
\Vert_{\mathcal{L}\left(H\right)} < 1$, when $k_0 \geq
\frac{1}{\theta} + 1$.

\hfill \hfill \textcolor{Aquamarine}{$\blacksquare $}
\end{proof}

\textcolor{Crimson}{\subsubsection{Proof of Lemma 2}}

\begin{proof}
Let us first assume that $x \geq y > 0$. 

From the Mean Value Theorem
applied over $e^z$, there exists $0 < \alpha < 1$ such that
$$\frac{e^{z+h} - e^z}{h} = e^{z + \alpha h}.$$ Taking $z = - x t$ and
$z + h = - yt$, we get the following inequalities:
\begin{equation}
\vert e^{- xt} - e^{-yt} \vert = \vert x - y \vert t e^{-x t +
\alpha \left(x - y \right)t} = \vert x - y \vert t e^{x t
\left(\alpha -1 \right)} e^{- y \alpha t}  \leq \vert x - y \vert t
e^{- y \alpha t} \leq\vert x - y \vert t. \nonumber %\label{A2:Supp:38}
\end{equation}

Similar inequalities are obtained for the case $y \geq x > 0,$ by
applying the Mean Value Theorem over the interval $\left[x,~y \right]$,
instead of $\left[y ,~x \right]$.

\hfill \hfill \textcolor{Aquamarine}{$\blacksquare$}
\end{proof}

\textcolor{Crimson}{\subsubsection{Proof of Lemma 3}}

\begin{proof}
 Considering  the indicator function $\boldsymbol{1}_{\cdot},$ it holds
\begin{eqnarray}
Y_n \vert Z_n \vert &=& Y_n \vert Z_n \vert \boldsymbol{1}_{\{\vert Z_n \vert
< a_n\}} + Y_n \vert Z_n \vert \boldsymbol{1}_{\{\vert Z_n \vert \geq a_n\}}
\leq Y_n a_n + Y_n \vert Z_n \vert \boldsymbol{1}_{\{\vert Z_n \vert \geq
a_n\}}, \label{A2:Supp:34}
\end{eqnarray}
\noindent where  $\left\lbrace a_n, \ n \in \mathbb{Z} \right\rbrace$ is a
sequence of positive numbers such that the event $\left\lbrace Y_n
\vert Z_n \vert \boldsymbol{1}_{\{\vert Z_n \vert \geq a_n\}} , \ n \in \mathbb{Z} \right\rbrace$
is equivalent to $\left\lbrace \vert Z_n \vert \geq a_n, \ n \in \mathbb{Z}
\right\rbrace$. From (\ref{A2:Supp:34}) and \textcolor{Crimson}{Lemma} \ref{A2:Supp:lempr1}, if we take
$a_n > \frac{\varepsilon}{2},$ for all $n \in \mathbb{Z},$ we get, for each $\varepsilon > 0$,
\begin{equation}
\mathcal{P} \left( Y_n \vert Z_n \vert \geq \varepsilon \right) 
\leq  \mathcal{P} \left( Y_n a_n \geq \frac{\varepsilon}{2} \right)
+ \mathcal{P} \left( \vert Z_n \vert \geq a_n  \right)  \leq
\mathcal{P} \left( Y_n a_n \geq \frac{\varepsilon}{2} \right) +
e^{-\theta a_{n}^{2}}. \label{A2:Supp:35}
\end{equation}

For $a_n = c \sqrt{\ln \left( n \right)} > \frac{\varepsilon}{2}$,
with $\frac{1}{\sqrt{\theta}} < c < + \infty$,
\begin{equation}
\displaystyle \sum_{n \in \mathbb{Z}} \mathcal{P} \left(\vert Z_n \vert \geq
a_n \right) \leq \displaystyle \sum_{n \in \mathbb{Z}} e^{-\theta a_{n}^{2}}
= \displaystyle \sum_{n \in \mathbb{Z}} \frac{1}{n^{\theta c^2}} < +\infty, \nonumber %\label{A2:Supp:36}
\end{equation}
\noindent which implies that $$\lim_{n\rightarrow
\infty}\mathcal{P}\left(\vert Z_n \vert \geq a_n \right)= 0$$  in
equation (\ref{A2:Supp:35}). 

On the other hand, since $\sqrt{\ln\left( n
\right)} Y_n \longrightarrow^{p} 0,$ for every $\varepsilon > 0,$
\begin{equation}
0=\lim_{n\rightarrow \infty}\mathcal{P} \left(\sqrt{\ln\left( n
\right)} Y_n\geq \frac{\varepsilon}{2} \right)= \lim_{n\rightarrow
\infty}\mathcal{P} \left( Y_n \frac{a_n}{c} \geq
\frac{\varepsilon}{2} \right). \nonumber %\label{A2:Supp:37}
\end{equation}
Thus, $Y_n \vert Z_n \vert \longrightarrow^{p} 0$.

\hfill \hfill \textcolor{Aquamarine}{$\blacksquare$}
\end{proof}

\textcolor{Crimson}{\section*{\textbf{Acknowledgments}}}

\textcolor{Aquamarine}{\textbf{This work has been supported in part by projects MTM2012-32674 and MTM2015--71839--P (co-funded by Feder funds), of the DGI, MINECO, Spain.}}

\vspace{0.5cm}
\renewcommand\bibname{\textcolor{Crimson}{\textit{\textbf{References}}}}

\bibliographystyle{dinat}
\bibliography{Biblio}

\end{document}